\documentclass[10pt]{article}
\usepackage[includeheadfoot,top=2 cm, bottom=2 cm, left=1.5 cm, right=1.5 cm]{geometry}
\usepackage{graphicx}
\usepackage{lipsum}
\usepackage{amsfonts,amsthm}
\usepackage{amssymb}
\usepackage{upgreek}
\usepackage{hyperref}
\usepackage[section]{algorithm}
\usepackage{algorithmicx,algcompatible,eqparbox}

\makeatletter
\renewcommand{\ALG@beginalgorithmic}{\small}
\makeatletter
\usepackage{multirow}
\usepackage{array}
\usepackage{mathtools}
\usepackage{color}
\usepackage{xcolor}
\usepackage{float}
\usepackage{subcaption}
\usepackage[noabbrev,capitalise]{cleveref}
\hypersetup{
	colorlinks,
	linkcolor={blue!},
	citecolor={blue!},
	urlcolor={blue!}
}

\usepackage[
sortcites,
backend=bibtex,
hyperref=true,
firstinits=true,
maxbibnames=200,
backref,
backrefstyle=none
]{biblatex} 
\renewbibmacro{in:}{}
\Crefname{proposition}{Proposition}{Propositions}
\crefname{equation}{}{}
\Crefname{equation}{}{}

\newcommand{\cond}{\textup{cond}}
\newcommand{\fl}{\textup{fl}}

\newcommand{\bDelta}{\mathbf{\Delta}}

\newcommand{{{\bxi}}}{\mathbf{\upxi}}
\newcommand{\Frob}{\mathrm{F}}

\newcommand{\ba}{\mathbf{a}}
\newcommand{\bp}{\mathbf{p}}
\newcommand{\bb}{\mathbf{b}}
\newcommand{\bu}{\mathbf{u}}
\newcommand{\be}{\mathbf{e}}
\newcommand{\bv}{\mathbf{v}}
\newcommand{\by}{\mathbf{y}}

\newcommand{\bq}{\mathbf{q}}
\newcommand{\bs}{\mathbf{s}}
\newcommand{\bz}{\mathbf{z}}
\newcommand{\bx}{\mathbf{x}}
\newcommand{\bw}{\mathbf{w}}

\newcommand{\br}{\mathbf{r}}
\newcommand{\bnull}{\mathbf{0}}
\newcommand{\bTheta}{\mathbf{\Theta}}

\newcommand{\bPhi}{\mathbf{\Phi}}
\newcommand{\bR}{\mathbf{R}}
\newcommand{\bA}{\mathbf{A}}

\newcommand{\bB}{\mathbf{B}}
\newcommand{\bhS}{\widehat{\mathbf{S}}}
\newcommand{\bhR}{\widehat{\mathbf{R}}}
\newcommand{\bhW}{\widehat{\mathbf{W}}}
\newcommand{\bhB}{\widehat{\mathbf{B}}}
\newcommand{\bhQ}{\widehat{\mathbf{Q}}}

\newcommand{\bhX}{\widehat{\mathbf{X}}}

\newcommand{\bhP}{\widehat{\mathbf{P}}}

\newcommand{\bhH}{\widehat{\mathbf{H}}}

\newcommand{\bhx}{\widehat{\mathbf{x}}}

\newcommand{\bhp}{\widehat{\mathbf{p}}}

\newcommand{\bX}{\mathbf{X}}
\newcommand{\bZ}{\mathbf{Z}}

\newcommand{\bH}{\mathbf{H}}

\newcommand{\bW}{\mathbf{W}}

\newcommand{\bS}{\mathbf{S}}
\newcommand{\bY}{\mathbf{Y}}
\newcommand{\bQ}{\mathbf{Q}}
\newcommand{\bC}{\mathbf{C}}
\newcommand{\bP}{\mathbf{P}}
\newcommand{\bI}{\mathbf{I}}

\newcommand{\bU}{\mathbf{U}}

\newcommand{\bV}{\mathbf{V}}

\newcommand{\bPi}{\mathbf{\Pi}}
\newcommand{\bLambda}{\mathbf{\Lambda}}
\newcommand{\bc}{\mathbf{c}}

\newtheorem{lemma}{lemma}[section]
\newtheorem{proposition}[lemma]{Proposition}
\newtheorem{corollary}[lemma]{Corollary}
\newtheorem{theorem}[lemma]{Theorem}
\newtheorem{remark}[lemma]{Remark}
\newtheorem{definition}[lemma]{Definition}
\newtheorem{assumptions}[lemma]{Assumptions}
\renewcommand{\thefootnote}{\fnsymbol{footnote}}
\numberwithin{equation}{section}

\begin{document}
\title{Randomized block Gram-Schmidt process for the solution of linear systems and eigenvalue problems.}
\author{Oleg Balabanov\footnotemark[1]~~and~~Laura Grigori\footnotemark[2]~\footnotemark[3]}	
\footnotetext[1]{{The work of this author was in part done at Alpines, Inria, Sorbonne Universit\'e, Universit\'e de Paris, CNRS, Laboratoire Jacques-Louis Lions, F-75012 Paris, France (olegbalabanov@gmail.com).}}
 \footnotetext[2]{
 	The work of this author was in part done at Alpines, Inria, Sorbonne Universit\'e, Universit\'e de Paris, CNRS, Laboratoire Jacques-Louis Lions, F-75012 Paris, France.}  
\footnotetext[3]{
Paul Scherrer Institute, Laboratory for Simulation and Modelling, 5232 PSI Villigen, and
\'Ecole Polytechnique F\'ed\'erale de Lausanne (EPFL), Institute of Mathematics, 1015 Lausanne, Switzerland (laura.grigori@epfl.ch).}

\renewcommand{\thefootnote}{\arabic{footnote}}

\date{}
\maketitle

\begin{abstract}
	This article introduces randomized block Gram-Schmidt process (RBGS) for QR decomposition. RBGS extends the single-vector randomized Gram-Schmidt (RGS) algorithm and inherits its key characteristics such as being more efficient and having at least as much stability as any deterministic (block) Gram-Schmidt algorithm. 

	Block algorithms offer superior performance as they are based on BLAS3 matrix-wise operations and reduce communication cost when executed in parallel. Notably, our low-synchronization variant of RBGS can be implemented in a parallel environment using only one global reduction operation between processors per block.		 
	Moreover, the block Gram-Schmidt orthogonalization is the key element in the block Arnoldi procedure for the construction of a Krylov basis, which in turn is used in GMRES, FOM and Rayleigh-Ritz methods for the solution of linear systems and clustered eigenvalue problems. In this article, we develop randomized versions of these methods, based on RBGS, and validate them on nontrivial numerical examples. 
	\end{abstract}


\section{Introduction} 

Let $\bW \in  \mathbb{R}^{n \times m}$ be a matrix with a moderately large number of columns, so that $m \ll n$.  We consider column-oriented block Gram-Schmidt (BGS) algorithms for {computing a} QR factorization of $\bW$:
$$\bW = \bQ \bR,$$ 
where $\bQ \in  \mathbb{R}^{n \times m}$ has $\ell_2$-orthonormal or very well-conditioned columns such that $\mathrm{range}(\bQ) = \mathrm{range}(\bW)$, and $\bR \in \mathbb{R}^{m \times m}$ is upper triangular with positive diagonal entries. Block algorithms are usually based on matrix-wise BLAS3 operations allowing proper exploitation of modern cache-based and high-performance computational architectures. The BGS orthogonalization forms a skeleton for block Krylov subspace methods for solving clustered eigenvalue problems as well as linear systems with multiple right-hand sides. It is also used in $s$-step, enlarged and other communication-avoiding Krylov subspace methods~\parencite{hoemmen2010communication,grigori2016enlarged}. Please see~\parencite{carson2020overview} and the references therein for an extensive overview of BGS variants, and~\parencite{zou2021gmres,saad2011numerical,swirydowicz2021low,stewart2002krylov,baker2006improving} for the underlying block Krylov methods.

{{The randomized Gram-Schmidt (RGS)} process for QR decomposition, based on the random sketching technique (see~\parencite{woodruff2014sketching,vershynin2018high,martinsson2020randomized} and the references therein), {was introduced in~\parencite{balabanov2021randomizedGS,balabanov2020randomized}.} It requires nearly half as many flops and data passes as the classical and modified Gram-Schmidt (CGS and MGS) processes. Moreover, in a parallel architecture, RGS maintains the number of synchronization points of CGS, and therefore reduces them by a factor of $\mathcal{O}(m)$ compared to MGS. Notably, RGS is significantly more stable than CGS and at least as stable as MGS, as it yields a well-conditioned Q factor even when $\cond(\bW) = \mathcal{O}(u^{-1})$, with $u$ representing the unit roundoff.} 
While the RGS algorithm can be already beneficial under unique precision computations, there are even more benefits that can be gained by working in two precisions: using a coarse unit roundoff for expensive high-dimensional operations and a fine unit roundoff elsewhere. In this case the stability of the RGS algorithm can be guaranteed for the coarse roundoff independent of the dominant dimension $n$ of the matrix. This property can be particularity useful for large-scale computations performed on low-precision arithmetic architectures. Another advantage of RGS is the ability of efficient a posteriori certification of {the factorization} without having to  estimate the condition number of a large-scale matrix, or to perform any other expensive large-scale computations.

{The essential feature of RGS algorithm is orthonormalization of a random sketch $\bTheta \bQ$ of $\bQ$ rather than the full factor, where $\bTheta \in \mathbb{R}^{k \times n}$ is a carefully chosen random matrix with typically $k = \mathcal{O}(m)$ rows that can be efficiently applied to vectors within the given architecture. The sketching matrix is designed to be an approximate isometry, or a $\varepsilon$-embedding, for $\mathrm{range}(\bQ)$ with high probability. The authors then applied RGS to Krylov methods for solving linear systems $\bA \bx = \bb$ and eigenvalue problems $\bA \bx = \lambda \bb$. They obtained randomized Arnoldi iteration providing a sketch-orthonormal Krylov basis $\bQ = [\bq_1, \hdots, \bq_{m}]$ that satisfies the Arnoldi identity  $\bA [\bq_1, \hdots, \bq_{m-1}] = \bQ \bH$, and the associated GMRES method providing a solution $\bx_{m-1}$ that minimizes the sketched residual $\| \bTheta (\bA \bx_{m-1} - \bb) \|_2$ over the Krylov subspace. The approximate orthogonality of $\bQ$ and the accuracy of $\bx_{m-1}$ compared to the classical GMRES solution directly follow from the fact that $\bTheta$ is a $\varepsilon$-embedding for the computed Krylov subspace.}

In this paper we propose a block {generalization} (RBGS) of the RGS process.
It has similar stability guarantees as its single-vector counterpart, and similar flops count {(nearly half of that of block CGS)}. At the same time, thanks to {the} block paradigm, it is better suited to modern computational architectures. In particular, the major operations in {the} RBGS algorithm can be implemented using cache-efficient BLAS3 subroutines and a reduced number of synchronizations between distributed processors. {Notably, we introduce a low-synchronization variant of RBGS, which stands out for its requirement of only one global synchronization between distributed processors per block.} {In addition, RBGS is only weakly sensitive to the accuracy of inter-block orthogonalization, and inherits the ability of efficient certification of {factorization}.} {Similarly to RGS, the RBGS has the advantage of performing dominant operations with a unit roundoff independent of the matrix's primary dimension.}

{Furthermore, we address the application of the RBGS algorithm to {block} Krylov methods. We introduce {a} block {generalization} of the randomized Arnoldi algorithm from~\parencite{balabanov2020randomized,balabanov2021randomizedGS}, which is subsequently used to develop {randomized block} GMRES, full orthogonalization method (FOM), and Rayleigh-Ritz (RR) method for solving linear systems and eigenvalue problems. In exact arithmetic, these methods can be interpreted as minimization of a sketched residual norm or imposing a sketched Galerkin orthogonality condition on the residuals. 
Furthermore, we also discuss an application of RBGS to Krylov $s$-step methods.}

{It is noticed that our RBGS algorithm can be augmented with a Cholesky QR step to provide a QR factorization with {a} Q factor that is not just well-conditioned but $\ell_2$-orthogonal to machine precision. {A detailed presentation of Cholesky QR and its properties can be found in~\parencite{fukaya2014choleskyqr2,yamamoto2015roundoff}. Such augmented procedure can be readily used to facilitate} the classical GMRES, {FOM} and RR approximations, or in any other applications.}

{In addition to~\parencite{balabanov2020randomized,balabanov2021randomizedGS} and the present article, there is another notable work by Nakatsukasa and Tropp~\parencite{nakatsukasa2021fast} on randomized Krylov methods,} {which was conducted independently and made available slightly earlier than the present work.}  {It is also worth noting that the sketched minimal-residual condition used in randomized GMRES  and the sketched Galerkin orthogonality condition used in randomized FOM and RR were proposed in~\parencite{balabanov2019randomized,balabanov2019randomized2},} {where they were applied to compute approximate solutions of linear parametric systems in a low-dimensional subspace.}

This article is organized as follows. The basic notations are explained in~\cref{prem}. We introduce a general BGS process and particularize it to {a} few classical variants in~\Cref{BGSprocess}. In~\cref{sketching} we present the basic ingredients of the random sketching technique and also extend the results from~\parencite{balabanov2021randomizedGS} concerning the effect of sketching on rounding errors in a matrix-matrix product. In~\cref{RBGS}, we present novel RBGS algorithms, {including the low synchronization variant which relies on randomized CholeskyQR, a matrix based version of a vector oriented algorithm presented in~\parencite{balabanov2021randomizedGS,balabanov2020randomized}. Randomized CholeskyQR is highly related to~\parencite{rokhlin2008fast}.}  
The stability of RBGS is analyzed in~\cref{stability}.  \Cref{applications} discusses the application of the methodology to solving clustered eigenvalue problems and linear systems, possibly with multiple right-hand sides. For this we develop the randomized block Arnoldi iteration and the associated GMRES, {FOM}, and RR methods. Nontrivial numerical experiments in~\cref{experiments} demonstrate great potential of our methodology. {The proofs of theorems and propositions from the stability analysis are deferred to~\cref{proofs}.} In~\cref{concl} we conclude the work and provide an avenue of future research.

{As supplementary material,  {we provide} a discussion on the efficient solution of sketched least-squares problems, which {is a component} of the RGS and RBGS algorithms, and provide an analysis of the accuracy of the randomized RR approximation from the paper.}

\subsection{Preliminaries} \label{prem}
Throughout this work, we use the following notations, which is an adaptation of the notations from~\parencite{balabanov2021randomizedGS} to block linear algebra. As in~\parencite{balabanov2021randomizedGS}, we denote vectors  $\bx$ by bold lowercase letters. A matrix composed of several column vectors $\bx_1, \hdots, \bx_k$ is denoted with the associated bold capital letter and a subscript specifying the number of columns, i.e., $\bX_k$. If $k$ is constant, this notation can be simplified to $\bX$. The $(i,j)$-th block of a block matrix $\bX$ is denoted by $\bX_{(i,j)}$ i.e, we have 
\begin{center} 
	$ \footnotesize \bX = \left [\begin{matrix} \bX_{(1,1)} & \bX_{(1,2)}& \cdots & \bX_{(1,p)} \\ \bX_{(2,1)} & \bX_{(2,2)} & \cdots & \bX_{(2,p)} \\ \vdots & \vdots & \ddots & \\ \bX_{(l,1)} & \bX_{(l,2)}  & & \bX_{(l,p)}\end{matrix} \right ], $
\end{center}
\normalsize
for some $p$ and $l$. When $l=1$ (i.e., matrix $\bX$ is partitioned column-wise), the notation $\bX_{(1,j)}$ can be simplified to $\bX_{(j)}$.
{The sub-block of $\bX$ composed of blocks $\bX_{(i,j)}$ with $N_1 \leq i \leq N_2$ and $M_1 \leq j \leq M_2$, is denoted by $\bX_{(N_1:N_2,M_1:M_2)}$.}
We denote the matrix $\bX_{(N_1:N_2,M_1:M_1)}$ by simply $\bX_{(N_1:N_2,M_1)}$. Moreover, if $l =1$, then matrix $\bX_{(1:1,M_1:M_2)}$ is denoted by $\bX_{(M_1:M_2)}$, and matrix $\bX_{(M_1:M_1)}$ by $\bX_{(M_1)}$. The minimal and the maximal singular values of $\bX$ are denoted by $\sigma_{min}(\bX)$ and $\sigma_{max}(\bX)$, and the condition number by $\cond(\bX)$. We let $\langle \cdot, \cdot \rangle$ and $\|\cdot\| = \sigma_{max}(\cdot)$ be the $\ell_2$-inner product and $\ell_2$-norm, respectively. $\| \cdot \|_\Frob$ denotes the Frobenius norm.  For two matrices (or vectors) $\bX$ and $\bY$, the relation  $ \bX \leq \bY $   indicates that the entries of $\bX$ satisfy $x_{i,j} \leq y_{i,j}$. Furthermore, for a matrix (or a vector) $\bX$, we denote by $|\bX|$ the matrix $\bY$ with entries $y_{i,j}= |x_{i,j}|$. We also let $\bX^\mathrm{T}$ and $\bX^{\dagger}$ respectively indicate the transpose and the Moore--Penrose inverse of $\bX$. Finally, we let $\bI$ be {the} identity matrix of size appropriate to the expression where this notation is used. 

A quantity or an arithmetic expression $X$ computed with finite precision arithmetic is denoted by $\fl(X)$ or $\hat{X}$. {In addition, we also use the following ``plus-minus'' notation. The $\pm$ symbol indicates that the left hand side is bounded from above by the right hand side with $\pm$ replaced by a plus, and from below with $\pm$ replaced by a minus. Moreover, $\mp$ symbol is used in conjunction with $\pm$.}

\subsection{BGS process} \label{BGSprocess}

Let matrix  $\bW \in \mathbb{R}^{n \times m}$ be partitioned
into $p$ blocks $\bW_{(i)} \in \mathbb{R}^{n \times m_p}$, with $1 \leq i \leq p$, $m = m_p p$:
\begin{center} $ \bW = \bW_{(1:p)}= \left [\begin{matrix} \bW_{(1)} \bW_{(2)} \ldots \bW_{{(p)}} \end{matrix} \right ]. $
\end{center}

{The} BGS process {for computing a QR factorization of $\bW$} proceeds recursively, at iteration $i$, selecting a new block matrix $\bW_{(i)}$ and orthogonalizing it with the previously orthogonalized blocks yielding a matrix $\bQ'_{(i)}$, followed by orthogonalization of $\bQ'_{(i)}$ itself.  This procedure is summarized in~\cref{alg:BGS}.

\begin{algorithm}[h] \caption{BGS process} \label{alg:BGS}
	\begin{algorithmic}
		\STATE{\textbf{Given:} $n \times m$ block matrix $\bW = \bW_{(1:p)}$, $m \leq n$ }
		\STATE{\textbf{Output}:  $n \times m$ factor $\bQ = \bQ_{(1:p)}$ and $m \times m$ upper triangular factor $\bR = \bR_{(1:p,1:p)}.$}
		\FOR{ $i = 1:p$} 
		\STATE{1. Compute a projection $\bQ'_{(i)} = \bPi^{(i-1)} \bW_{(i)}$  (also yielding $\bR_{(1:i-1,i)}$)}.
		\STATE{2. Compute QR factorization $\bQ_{(i)}\bR_{(i,i)} = \bQ'_{(i)}$ with suitable efficient routine.}
		\ENDFOR
	\end{algorithmic}
\end{algorithm}

For standard methods, the projector $\bPi^{(i-1)}$ in~\cref{alg:BGS} is taken as an approximation to $\ell_2$-orthogonal projector onto $\mathrm{range}(\bQ_{(1:i-1)})^{\perp}$. For the  classical BGS process (BCGS), one chooses $\bPi^{(i-1)}$ as $$\bPi^{(i-1)} = \bI -  \bQ_{(1:i-1)} \bQ_{(1:i-1)}^\mathrm{T}.$$ Whereas, for the modified BGS (BMGS) process, we have
\begin{equation} \label{eq:blockMGS}
	\bPi^{(i-1)} = (\bI -  \bQ_{(i-1)} (\bQ_{(i-1)})^\mathrm{T}) (\bI -  \bQ_{(i-2)} (\bQ_{(i-2)})^\mathrm{T}) \ldots (\bI -  \bQ_{(1)} (\bQ_{(1)})^\mathrm{T}).
\end{equation}
In infinite precision arithmetic, these two projectors are equivalent, and are exactly equal to the  $\ell_2$-orthogonal projector onto $\mathrm{range}(\bQ_{(1:i-1)})^{\perp}$, since, by construction, $\bQ_{(1:i-1)}$ is an orthogonal matrix.  However, in the presence of rounding errors these projectors may cause instability. In the first case, the condition number of $\bQ$ can grow as $\cond(\bW)^{2}$ or even worse and thus requires special treatment~\parencite{carson2021stability}. While in the second case, $\cond(\bQ)$ can grow as $\cond(\bW) \max_{1\leq j \leq p}\cond(\bW_{(j)})$ unless step 2 is unconditionally stable~\parencite{barlow2019block,carson2020overview}. The stability of these processes can be improved by re-orthogonalization, i.e., by running the inner loop twice. In particular, it can be shown that the BCGS process with re-orthogonalization, here denoted by BCGS2, yields an almost orthogonal Q factor as long as the matrix $\bW$ is numerically {full rank}~\parencite{barlow2013reorthogonalized}. 

Inter-block orthogonalization is an important step in BGS algorithm. For standard versions, this step can be performed with any suitable efficient and stable routine for ($\ell_2$-)QR factorization of tall and skinny matrices, such as TSQR from~\parencite{demmel2012communication}, a Gram-Schmidt QR applied $l$ times, Cholesky QR applied $l$ times, and others.  Typically, such QR factorization takes only a fraction of the overall computational cost.

\subsection{Random sketching} \label{sketching}
{In this subsection we recall the basic notions of the random sketching technique.} Let $\bTheta \in {\mathbb{R}}^{k\times n}$ with $k \leq n$ be a sketching matrix such that the  associated sketched product $\langle \bTheta \cdot, \bTheta \cdot \rangle $ approximates well the $\ell_2$-inner product between any two vectors in the subspace (or subspaces) of interest $V \subset \mathbb{R}^{n}$. Such matrix $\bTheta$ is referred to as {a} $\varepsilon$-embedding for $V$, as defined below.

\begin{definition}
	For $\varepsilon<1$, we say that the sketching matrix $\bTheta$ is {a} $\varepsilon$-embedding for subspace $V$ {(or matrix $\bV$ spanning $V$)}, if it satisfies the following relation:
	\begin{equation} \label{eq:isometry}
	\forall \bx,\by \in V,~~ | \langle \bx, \by \rangle - \langle \bTheta \bx, \bTheta \by \rangle | \leq \varepsilon \| \bx \| \|\by\|.  
	\end{equation} 
\end{definition}

Furthermore, we assume that $\bTheta$ is chosen at random from a certain distribution, such that it satisfies~\cref{eq:isometry} for any fixed $d$-dimensional subspace with high probability (see~\cref{def:oblemb}).
\begin{definition} \label{def:oblemb}
	A random sketching matrix $\bTheta$ is called a $(\varepsilon,\delta,d)$ oblivious $\ell_2$-subspace embedding, if for any fixed $V \subset \mathbb{R}^n$ of dimension $d$, it satisfies 
	$$\mathbb{P}(\text{$\bTheta$ is a $\varepsilon$-embedding for $V$}) \geq 1-\delta. $$
\end{definition}

\begin{corollary}\label{prop:Cn2sqrtn}
	If $\bTheta \in \mathbb{R}^{k \times n}$ is a $(\varepsilon, \delta/n, 1)$ oblivious $\ell_2$-subspace embedding, then with probability at least $1-\delta$, we have
	$$\|\bTheta\|_\Frob \leq\sqrt{(1+\varepsilon)n}. $$
\end{corollary}

\begin{corollary}\label{prop:skcond} If $\bTheta \in \mathbb{R}^{k \times n}$ is {an} $\varepsilon$-embedding for $\bV$, then the singular values of $\bV$ are bounded by $$ (1+\varepsilon)^{-1/2} \sigma_{min}(\bTheta\bV)  \leq \sigma_{min}(\bV) \leq \sigma_{max}(\bV) \leq  (1-\varepsilon)^{-1/2} \sigma_{max}(\bTheta\bV).$$
	\vspace*{-0.5cm}	
\end{corollary}

{The proofs for~\cref{prop:Cn2sqrtn,prop:skcond} can be found for instance in~\parencite{balabanov2021randomizedGS}.}

In recent years, several distributions of $\bTheta$ have been proposed that satisfy~\cref{def:oblemb} and have a small first dimension $k$ {that depends at most logarithmically} on the dimension $n$ and probability of failure $\delta$. Among them, the most suitable distribution should be selected depending
on the problem and computational architecture.
The potential of random sketching is here realized on (rescaled) Rademacher matrices and Subsampled Randomized Hadamard Transform (SRHT). The entries of a Rademacher matrix are i.i.d. random variables satisfying $\mathbb{P}(\theta_{i,j}  = 1/\sqrt{k})=\mathbb{P}(\theta_{i,j}  = - 1/\sqrt{k})=1/2$.  Rademacher matrices can be efficiently multiplied by vectors and matrices through the proper exploitation of computational resources such as cache or distributed machines.
For $n$, which is a power of $2$, SRHT is defined as a product of {a} diagonal matrix of random signs with {a} Walsh-Hadamard matrix, followed by {an} uniform sub-sampling matrix {scaled by} $1/\sqrt{k}$. For a general $n$, the SRHT has to be combined with zero padding to make the dimension a power of $2$. Random sketching with SRHT can reduce  the complexity of an algorithm. Products of SRHT matrices with vectors require only $n \log_2 n$ flops  using  the  fast  Walsh-Hadamard  transform or $2n\log_2(k+ 1)$ flops following the methods  in~\parencite{ailon2009fast}.   Furthermore, for both distributions, the usage of a seeded random number generator can allow efficient {storage and application} of $\bTheta$. It follows {from~\parencite{woodruff2014sketching,tropp2011improved}} that the rescaled Rademacher distribution, and SRHT (possibly with zero padding) respectively are $(\varepsilon, \delta, d)$ oblivious $\ell_2$-subspace embeddings, if they have a sufficiently large first dimension: if {$k \geq 7.87 \varepsilon^{-2}( {6.9} d + \log (1/\delta))$ for Rademacher matrices, or if $k \geq 2( \varepsilon^{2} - \varepsilon^3/3)^{-1} \left (\sqrt{d}+ \sqrt{8 \log(6 n/\delta)} \right )^2 \log (3 d/\delta)$ for SRHT. {A complete proof of this fact can be found for instance in \parencite{balabanov2019randomized}.}}

\subsection{Effect of random sketching on rounding errors} \label{sketchrounding}
Let us now discuss {an important} result from~\parencite[Section 2.2]{balabanov2021randomizedGS} characterizing the rounding errors in a sketched matrix-vector product, which can be easily extended to matrix-matrix products. In short, this result states that multiplying $\bTheta$ by a matrix-vector product $\bhx = \mathrm{fl}(\bY \bz)$ does not increase the rounding error bound of $\bhx$ by more than a small factor. This property is {essentially} a sketched version of the standard ``rule of thumb'' {of rounding analysis~\parencite{higham2002accuracy}}. It was rigorously proven in~\parencite{balabanov2021randomizedGS} for the probabilistic rounding model, \parencite[Model 4.7]{connolly2020stochastic}. The extension of this result to matrix-matrix products is provided below. 

Let us fix a realization of an oblivious $\ell_2$-subspace embedding  $\bTheta \in \mathbb{R}^{k \times n}$ of sufficiently large size, and consider a matrix-matrix product $\bX = \bY \bZ$, with $\bY \in \mathbb{R}^{n\times m},~\bZ \in \mathbb{R}^{m \times l} $,
computed in finite precision arithmetic with unit roundoff $u <0.01/m$.  The following results can be derived directly from~\parencite[Section 2.2]{balabanov2021randomizedGS} with the observation that each column of $\bX$ is a matrix-vector product: $\bx_i = \bY \bz_i,~1 \leq i \leq l.$

We have,
\begin{equation} \label{eq:roundoffYz}
	|\bX - \bhX| \leq  \bU,
\end{equation}
for some matrix $\bU$ describing the ``worst-case scenario'' rounding error.
In general, the standard analysis (e.g., see~\parencite{higham2002accuracy}) gives the bound
\begin{equation} 
	|\bX - \bhX| \leq  \frac{ m u}{1-mu} |\bY| |\bZ| \leq 1.02 m u |\bY| |\bZ|,
\end{equation}
which means that  $\bU= 1.02 m u |\bY| |\bZ|$ satisfies~\cref{eq:roundoffYz}. Furthermore, as in~\parencite{balabanov2021randomizedGS}, if $\bY \bZ$ represents a sum of a matrix-matrix product with a matrix, i.e., $\bY \bZ = \bY' \bZ'+\bH$, then one can take
$$\bU =  1.02 u (|\bH| + m |\bY'| |\bZ'|).$$

Let us now address bounding the norm of the rounding error after sketching $\bhX$. We have the following ``worst-case scenario'' bound
\begin{equation} 
	\|\bTheta(\bX - \bhX) \|_\Frob \leq \|\bTheta \|  \|\bX - \bhX \|_\Frob \leq \|\bTheta \|  \|\bU \|_\Frob,
\end{equation}
which, combined with~\cref{prop:Cn2sqrtn}, implies that if $\bTheta$ is $(\varepsilon, \delta/n, 1)$ oblivious $\ell_2$-subspace embedding, then $\| \bTheta(\bX - \bhX) \|_\Frob \leq \sqrt{1+\varepsilon} \sqrt{n} \|\bU \|_\Frob$ holds with probability at least $1-\delta$. The next important point is that, as was argued in~\parencite{balabanov2021randomizedGS}, this bound is pessimistic and can be improved by a factor of $\mathcal{O}{(\sqrt{n})}$ by exploiting statistical properties of rounding errors. {For instance one can rely on the assumption that the rounding errors in elementary arithmetic operations are mean-independent random variables~\parencite{connolly2020stochastic}.}  We summarize the results from~\parencite[Theorem 2.5 and Corollary 2.6]{balabanov2021randomizedGS} below. 

\begin{corollary} \label{thm:thetadeltax2}
	{Consider a matrix-matrix product $\bX = \bY \bZ$, with $\bY \in \mathbb{R}^{n\times m},~\bZ \in \mathbb{R}^{m \times l} $, computed under {probabilistic} rounding model, where the rounding errors due to elementary arithmetic operations are mean-independent random variables with zero mean.}	
	Furthermore assume that the errors are bounded so that it holds,
	$$|\bX - \bhX| \leq \bU,$$
	{where $\bU$ is a deterministic matrix representing the worst-case scenario rounding error.}
	If $\bTheta$ is a $(\varepsilon/4,l^{-1}\binom{n}{d}^{-1}\delta, d)$ oblivious $\ell_2$-subspace embedding, with $d = 4.2 c^{-1} \log (4/\delta)$, where $c \leq 1$ is some universal constant, then 
	\begin{equation}\label{eq:thetadeltax2}
		\| \bTheta(\bX - \bhX) \|_\Frob \leq \sqrt{1+\varepsilon} \|\bU \|_\Frob
	\end{equation}
	holds with probability at least $1-2\delta$. 
\end{corollary}
By using {the bounds from~\cref{sketching}} we deduce that for $l \leq n$ the relation~\cref{eq:thetadeltax2} is satisfied with probability at least $1-2\delta$ if $\bTheta$ is a Rademacher matrix with {$\mathcal{O}(\log(n)\log(1/\delta))$} rows or SRHT matrix with $\mathcal{O}(\log^2(n)\log^2(1/\delta))$ rows. 

{The fact that the bound~\cref{eq:thetadeltax2} is independent of (the high) dimension $n$ implies that, in practice,  {products of matrices} in randomized algorithms can be performed with a unit roundoff that is independent of $n$.}

\section{RBGS process} \label{RBGS}
Consider BGS algorithms with projectors $\bPi^{(i-1)}$ {that respect the relation}
$$ { \bPi^{(i-1)} \bW_{(i)} = \bW_{(i)} - \bQ_{(1:i-1)} \bX,} $$
where $\bX = \bR_{(1:i-1,i)}$ is computed from $\bW^{(i)}$ and $\bQ_{(1:i-1)}$. Standard algorithms take $\bX$ as an approximate solution to the following minimization problem:
\begin{equation} \label{eq:RBGSminimization}
\min_{\bY} \left \| \bQ_{(1:i-1)}   \bY  - \bW_{(i)} \right \|_\Frob.
\end{equation}
For instance, BCGS algorithm approximates the exact solution to~\cref{eq:RBGSminimization} by $\bQ_{(1:i-1)}^\mathrm{T}\bW_{(i)}$, whereas BMGS projector~\cref{eq:blockMGS} improves this approximation under finite precision arithmetic, though it may cause a computational overhead in terms of inter-processor communication and operation with cache/RAM.

As proposed in~\parencite{balabanov2021randomizedGS}, a reduction of computational cost and/or improvement of {the} stability of {the} Gram-Schmidt process can be obtained with a projector that gives a Q factor orthonormal with respect to  $\langle \bTheta \cdot, \bTheta \cdot \rangle $ instead of the $\ell_2$-inner product as in standard methods. In our case, this corresponds to taking $\bX$ as an approximate solution to 
\begin{equation} \label{eq:skRBGSminimization}
	\min_{\bY} \left \| \bTheta \bQ_{(1:i-1)}\bY  - \bTheta\bW_{(i)} \right \|_\Frob,
\end{equation}
which is a $k$-dimensional block least-squares problem. Since the sketching dimension satisfies $k \ll n$, a very accurate solution to~\cref{eq:skRBGSminimization} should be more efficient to compute than even the {cheapest (i.e, BCGS)} solution to~\cref{eq:RBGSminimization}. Furthermore, with the right choice of random sketching matrices, the precomputation of the sketches of $\bQ_{(1:i-1)}$ and $\bW_{(i)}$ should also have only a minor cost compared to that of standard high-dimensional operations. If~\cref{eq:skRBGSminimization} is very small, its solution can be obtained with a direct solver based, for example, on {Householder or Givens QR factorization}, requiring a cubic complexity. {If the problem has a moderate size, so that each direct solution has a considerable computational cost, it can be beneficial to recycle the QR factorization of $\bTheta \bQ_{(1:i-2)}$ computed at iteration $i-1$ to get the QR factorization of  $\bTheta \bQ_{(1:i-1)}$. Alternatively, we may use the fact that} {the} {matrix  $\bTheta\bQ_{(1:i-1)}$ is almost orthonormal, which implies applicability of iterative solvers running in quadratic complexity. Several such solvers are discussed in~\cref{lssol}.}

To produce a QR factorization of $\bW$ with respect to $\langle \bTheta \cdot, \bTheta \cdot \rangle $,  the inter-block orthogonalization of $\bQ'_{(i)}$ in step 2 of~\cref{alg:BGS} has to provide a Q factor $\bQ_{(i)}$ orthonormal with respect to $\langle \bTheta \cdot, \bTheta \cdot \rangle $. An efficient procedure for this task can be based on {sketched Cholesky QR  factorization} (RCholeskyQR), which consists in obtaining the R factor $\bR_{(i)}$ by a regular QR of the sketch $\bTheta \bQ'_{(i)}$, and then retrieving $\bQ_{(i)}$ with forward substitution: $\bQ_{(i)} = \bQ'_{(i)} \bR_{(i)}^{-1}$. In this case, it can be beneficial to obtain the sketch of $\bQ'_{(i)}$ from  $ \bTheta \bQ'_{(i)} = \bTheta \bW_{(i)} - \bTheta \bQ_{(1:i-1)} \bR_{(1:i-1,i)}$ rather than by multiplying $\bTheta$ with $\bQ'_{(i)}$, as this saves a global synchronization between processors. Another option would be to simply perform the inter-block orthogonalization with a single-vector RGS algorithm. See~\cref{interQR} for more details.

The RBGS process using the random sketching projector is depicted in~\cref{alg:RBGS}.  
\begin{algorithm}[h] \caption{RBGS algorithm (RBGS)} \label{alg:RBGS}
	\begin{algorithmic}
		\STATE{\textbf{Given:} $n \times m$  block matrix $\bW = \bW_{(1:p)}$, and $k \times n$ matrix $\bTheta$, $m \leq k \ll n$.} 
		\STATE{\textbf{Output}:  $n \times m$ factor $\bQ = \bQ_{(1:p)}$ and $m \times m$ upper triangular factor $\bR = \bR_{(1:p,1:p)}.$}
		\FOR{$i = 1:p$} 
		\STATE{1. Sketch $\bW_{(i)}$: $\bP_{(i)} = \bTheta \bW_{(i)}$.} \COMMENT{macheps: $u_{fine}$}
		\STATE{2. Solve small block least-squares problem:
			\\~~~~~~~~~~~~~~~~~~~~~~~~~~~~~~~{$\bR_{(1:i-1,i)} = \arg \min_{\bY} \left \| \bS_{(1:i-1)} \bY  - \bP_{(i)} \right \|_\Frob.$} \normalsize } \COMMENT{macheps: $u_{fine}$} 
		\STATE{3. Compute projection of $\bW_{(i)}$:  $\bQ'_{(i)} = \bW_{(i)} - \bQ_{(1:i-1)} \bR_{(1:i-1,i)}.$} \COMMENT{macheps: $u_{crs}$}
		\STATE{4-5. Compute QR fact. $\bQ_{(i)}\bR_{(i,i)} = \bQ'_{(i)}$  with respect to $\langle \bTheta \cdot, \bTheta \cdot \rangle$}. \COMMENT{macheps: $u_{fine}$}
		\STATE{Compute sketch of $\bQ_{(i)}$: $\bS_{(i)} = \bTheta \bQ_{(i)}$.}  \COMMENT{macheps: $u_{fine}$}
		\ENDFOR
		\STATE{6. (Optional) compute $\Delta^{(p)}  = \|\bI - \bS^\mathrm{T} \bS \|_\Frob$ and $\tilde{\Delta}^{(p)}= \frac{\|\bP - \bS\bR \|_\Frob}{\|\bP\|_\Frob}$. } \COMMENT{macheps: $u_{fine}$}	
	\end{algorithmic}
\end{algorithm}

{At the iteration $i = 1$ in~\cref{alg:RBGS} we used the notation that $\bR_{(1:i-1,i)}$ is a $0 \times m_p$ matrix and $\bQ_{(1: i-1)} $ is a $ n \times 0$ matrix, so that $\bQ'_{(i)} = \bW_{(i)}$.}  \Cref {alg:RBGS} is presented under multi-precision arithmetic using two unit roundoffs, like its single-vector counterpart from~\parencite{balabanov2021randomizedGS}. The working precision is represented by a coarse roundoff $u_{crs}$. It is used for standard high-dimensional operations in step 3, which determine the overall computational cost.   All other {(inexpensive)} operations in~\cref {alg:RBGS} are computed with a fine unit roundoff $u_{fine}$, $ u_ {fine} \leq u_{crs}$. It is shown in~\cref{stability} that \cref{alg:RBGS} is stable if $u_{crs} \leq \mathcal{O}(\cond{(\bW)} m^{-2})$, which is a very mild condition on $u_{crs}$. The fact that this bound is independent of the high dimension $n$ explains potential of our methodology for large-scale problems computed on low-precision arithmetic architectures. Furthermore, according to our numerical experiments, the RBGS algorithm can be sufficiently stable even when $u_{crs}$ is larger than $\mathcal{O}(\cond{(\bW)} m^{-2})$. {In such cases, the stability of the algorithm can be certified by a posteriori bounds given by the quantities $\Delta^{(p)}$ and $\tilde{\Delta}^{(p)},$ computed in (optional) step 6.} {These bounds provide stability certification if $u_{crs} = \mathcal{O}(m^{-3/2})$, which is a milder condition {than the one for a priori guarantees}, and is  independent not only of $n$ but also {of} $\cond{(\bW)}$.}

{The stability guarantees of \cref{alg:RBGS} executed in unique precision can be obtained directly from the analysis of the multi-precision algorithm by taking $u_{fine} = u$ and $u_{crs} = F(m, n)u$, where $F(m,n)$ is some polynomial of low degree.}

In terms of performance, the SRHT-based RBGS algorithm requires about half the flops and data passes of the cheapest standard BGS algorithm, which is BCGS. Its computational cost is defined by $p$ well-parallelizable BLAS3 operations. 
{Furthermore, with a suitable choice of inter-block QR factorization in steps 4-5, RBGS requires only one global reduction operation per block (see~\cref{interQR1} for details). This version of RBGS is particularly noteworthy in the realm of ``one-synchronization'' block Gram-Schmidt algorithms~\parencite{carson2022block,lund2022adaptively} in particular due to its provable stability characteristics.}

\begin{remark} \label{postproc}
	If necessary, the output of the RBGS algorithm can be  post-processed with Cholesky QR to provide a QR factorization  with {a} Q factor, which is not only well-conditioned but is $\ell_2$-orthonormal up to machine precision.  More specifically, we can compute an upper-triangular matrix $\bR'$ such that ${\bR'}^\mathrm{T} \bR' =  \bQ^\mathrm{T} \bQ$ with {a} Cholesky decomposition, and consider    
	$$ \bQ\leftarrow \bQ \bR'^{-1} \text{ and }  \bR \leftarrow \bR' \bR . $$
	{The computational cost of this procedure} is dominated by the computation of  $\bQ^\mathrm{T} \bQ$ and, possibly, $\bQ \bR'^{-1}$. In terms of flops, the computational cost of computing $\bQ^\mathrm{T} \bQ$ is similar to the cost of the entire RBGS algorithm. Though as a single BLAS3 operation, it is better suited for cache-based and  parallel computing. 
	The computation of $\bQ \bR'^{-1}$ can be omitted when the application allows operating with the Q factor in an implicit form, {for instance, in {the} Arnoldi iteration or RR algorithm.} Otherwise, this product can be computed by well-parallelizable forward substitution, or even {by} direct inversion, allowing a second BLAS3 multiplication. {The stability guarantees of the QR factorization obtained after the Cholesky QR step follow directly from the fact that the RBGS algorithm produces a well-conditioned Q factor along with the standard numerical stability guarantees of Cholesky QR.}
\end{remark}

\subsection{Solution of block least-squares problem in step 2} \label{lssol}
As already pointed out, the stability of~\cref{alg:RBGS} strongly depends on the stability of the least-squares solver used in step 2. In particular, in our analysis in~\cref{stability}  we require the solution $\bhX=\bhR_{(1:i-1,i)}$ to satisfy the following backward-stability condition.

\begin{assumptions}[Backward-stability of step 2 of RBGS] \label{thm:asssolver}
	For each $1\leq j \leq m_p$, the $j$-th column $\bhx$ of $\bhX$ and the $j$-th column $\bhp$ of $\bhP_{(i)}$ satisfy
	\begin{equation}
	\bhx= \arg \min_{\by} \left \| (\bhS_{(1:i-1)}+\bDelta \bS) \by  - (\bhp + \bDelta \bp  ) \right \|,
	\end{equation}
	where the perturbations $\bDelta \bS$ and $\bDelta \bp$ are such that
	$$  \| \bDelta \bS\|_\Frob \leq 0.01 u_{crs}   \|  \bhS_{(1:i-1)} \|,~~   \| \bDelta \bp \| \leq 0.01 u_{crs}   \|\bhp\|.$$ 
	Note that $\bDelta \bS$ and $\bDelta \bp$ may depend on $(i,j)$.
\end{assumptions} 

The condition in~\cref{thm:asssolver}  can be met by standard direct solvers based on Householder transformation or Givens rotations and a sufficiently large gap between $u_{crs}$ and $u_{fine}$, which follows from~\parencite[Theorems 8.5, 19.10 and 20.3]{higham2002accuracy} and their proofs. However, direct solvers require a cubic complexity and can become too expensive even for relatively small $m$ and $k$. {A remedy can be to re-use the QR factorization of $\bS_{(1:i-2)}$ computed at iteration $i-1$ to get the QR factorization of  $\bS_{(1:i-1)}$.}
{Another option is to} appeal to \emph{iterative} methods that can exploit the approximate orthogonality of $\bS_{(1:i-1)}$ to speedup the computations {such as the Richardson iterations $\bX \leftarrow \bX + \bS_{(1:i-1)}^\mathrm{T} \left (\bP_{(i)} -\bS_{(1:i-1)} \bX \right )$ (that can be  viewed as CGS or BCGS reorthogonalizations), MGS or BMGS reorthogonalizations, Conjugate Gradient or GMRES methods applied to the normal system of equations   
	$\left (\bS_{(1:i-1)}^\mathrm{T} \bS_{(1:i-1)} \right ) \bX =  \bS_{(1:i-1)}^\mathrm{T} \bP_{(i)}$. 
	A discussion of these methods can be found in the supplement to the article.}

\subsection{Inter-block QR factorization in steps 4--5} \label{interQR}
Let us now provide four ways of efficient and stable inter-block orthogonalization in steps 4--5 of the RBGS algorithm. 

\subsubsection{Single-vector RGS algorithm}
First, the inter-block orthogonalization can be readily performed with the single-vector RGS algorithm from~\parencite{balabanov2021randomizedGS} computed with unique roundoff $u = u_{fine} = F(m,n) u_{crs}$, where $F(m,n)$ is a low-degree polynomial.  This algorithm provides a QR factorization that satisfies~\cref{eq:sRGS1}, which follows directly by construction. Furthermore, it should also satisfy~\cref{eq:sRGS2}, since $\|\bTheta\|_\Frob \leq \sqrt{2n}$ holds with high probability (according to~\cref{prop:Cn2sqrtn}).   
\begin{subequations} \label{eq:sRGS}
	\begin{align}
		&|{\bhQ}'_{(i)} - {\bhQ}_{(i)} \bhR_{(i,i)}|   \leq 0.1 u_{crs} |{\bhQ}_{(i)}| | \bhR_{(i,i)} | \label{eq:sRGS1} \\
		&\|\bTheta({\bhQ}'_{(i)} - {\bhQ}_{(i)} \bhR_{(i,i)}) \|_\Frob  \leq 0.1 u_{crs} \| {\bhQ}_{(i)} \| \| \bhR_{(i,i)}\| \label{eq:sRGS2}
	\end{align}
\end{subequations}
As shown in~\parencite{balabanov2021randomizedGS}, if $\bTheta$ is $\varepsilon$-embedding for ${\bhQ}'_{(i)}$, then {the} RGS algorithm satisfies
\begin{equation} \label{eq:intercond}
1-0.1 u_{crs} \cond(\bhQ'_{(i)}) \leq  \sigma_{min}(\bTheta \bhQ_{(i)}) \leq \sigma_{max}(\bTheta \bhQ_{(i)}) \leq 1+ 0.1 u_{crs} \cond(\bhQ'_{(i)}).
\end{equation}
Moreover,  in this case $\bTheta$ is guaranteed to be a $\varepsilon'$-embedding for ${\bhQ}_{(i)}$, with $\varepsilon' = 2\varepsilon + u_{crs}$. Though, in practice, this property holds for smaller values $\varepsilon'$, say, $\varepsilon + u_{crs}$.

According to~\parencite{balabanov2021randomizedGS}, the single-vector RGS algorithm applied to an inter-block has the cost of $2 m_p^2 n$ flops and $m_p n$ memory {units}, which can be considered negligible compared to the complexity and memory consumption of other computations such as the standard high-dimensional operations in~\cref{alg:RBGS}. The RGS algorithm, however, requires $m_p$ global synchronizations between distributed processors, which can dominate the computational costs in parallel architectures. In such cases, it is necessary to appeal to {other} approaches for inter-block QR factorization, some of which are described below. 

\subsubsection{RCholeskyQR factorization} 

Another way to perform inter-block QR factorization of matrix $\bQ'_{(i)}$ with respect to the sketched inner product is to appeal to the following sketched version of Cholesky QR{, called RCholeskyQR.} We can first compute the R factor $\bR_{(i,i)}$ by performing a ($\ell_2$-)QR factorization of a small matrix $\bS'_{(i)} = \bTheta \bQ'_{(i)}$ with any suitable stable routine. Then the Q factor is retrieved by computing $\bQ_{(i)} = \bQ'_{(i)} ({\bR_{(i,i)}})^{-1}$ with forward substitution (see~\cref{alg:CholQR}). {We note that the vector oriented version of RCholeskyQR factorization, where a QR of $\bS'_{(i)}$ and the forward substitution $\bQ'_{(i)} ({\bR_{(i,i)}})^{-1}$ are performed column by column, was introduced first in~\parencite[Remark 2.10]{balabanov2020randomized,balabanov2021randomizedGS}. For completeness, we present this version in \Cref{alg:vCholeskyQR}, particularly useful when the vectors become available one at a time. As notation, the sub-block $\bX[N_1:N_2, M_1:M2]$ is formed by the elements in rows $N_1$ to $N_2$ and columns $M_1$ to $M_2$ of $\bX$, while the sub-block $\bX[:, M_1:M2]$ is formed by all the rows and columns $M_1$ to $M_2$ of $\bX$. We refer to this algorithm as vRCholeskyQR. The difference between vRCholeskyQR and RGS lies in step 4 of~\Cref{alg:vCholeskyQR}, where in RGS the sketch is computed directly from the columns of $\bQ_{(i)}$ before scaling, as $\widetilde{\bs}'_j = \bTheta \widetilde{\bq}'_j$. }

{In exact arithmetic RCholeskyQR provides the same output as {the} RGS algorithm from~\parencite{balabanov2020randomized,balabanov2021randomizedGS}. This in particular implies that the Q factor $\bQ_{(i)}$ obtained by RCholeskyQR is very well conditioned with a high probability.
Moreover, it was shown more recently in~\parencite{balabanov2022randomized} that RCholeskyQR has finite-precision stability guarantees similar to those of RGS given by~\cref{eq:sRGS,eq:intercond}.  However, in practice RCholeskyQR may provide less numerical stability than RGS, as is indicated in~\parencite[Remark 2.10]{balabanov2020randomized,balabanov2021randomizedGS}. }

{RCholeskyQR algorithm has a high relation to the preconditioning technique for overdetermined least-squares problems developed in~\parencite{rokhlin2008fast}. The difference is that in~\parencite{rokhlin2008fast}, $\bQ'_{(i)} ({\bR_{(i,i)}})^{-1}$ is not computed explicitly, but rather operated with as a function providing products with vectors, and that ${\bR_{(i,i)}}$ is not computed using a regular QR of $\bS'_{(i)}$ but rather a QR with column-pivoting.} 
{RCholeskyQR is also considered in \parencite{nakatsukasa2021fast} to approximately orthogonalize a basis, where it is referred to as ``basis whitening''.}

\begin{algorithm}[h] \caption{Vector oriented RCholeskyQR (vRCholeskyQR)} \label{alg:vCholeskyQR}
	\begin{algorithmic}
		\STATE{\textbf{Given:} $\bQ'_{(i)}$, $\bTheta$.} 
		\STATE{\textbf{Output}: QR fact. $\bQ'_{(i)} = \bQ_{(i)}\bR_{(i,i)}$, {where $\bR_{(i,i)}$ is upper triangular} and $\bQ_{(i)}$ is orthonormal with respect to $\langle \bTheta \cdot, \bTheta \cdot \rangle$}; a sketch $\bS_{(i)} = \bTheta \bQ_{(i)}$. For clarity, let $\bV := \bQ'_{(i)} = \begin{bmatrix} \bv_1, \ldots, \bv_{m_p} \end{bmatrix}$, $\widetilde{\bQ}:= \bQ_{(i)} = \begin{bmatrix} \widetilde{\bq}_1, \ldots, \widetilde{\bq}_{m_p} \end{bmatrix}$, $\widetilde{\bR} := \bR_{(i,i)}$, and $\widetilde{\bS} := \bS_{(i)} = \begin{bmatrix} \widetilde{\bs}_1, \ldots, \widetilde{\bs}_{m_p} \end{bmatrix}$	
		\FOR{$j = 1:m_p$} 
		\STATE{1. Sketch $\bv_j$: $\bp_j = \bTheta \bv_j$.} \COMMENT{macheps: $u_{fine}$}
		\STATE{2. Solve least-squares problem: \\~~~~~~~~~~{$\widetilde{\bR}{[1:j-1,j]} =  \arg \min_\by \| \widetilde{\bS}{[:,1:j-1]} \cdot \by - \bp_j \|.$} }
		\COMMENT{macheps: $u_{fine}$}
		\STATE{3. Compute projection of ${\bv}_j$: \\~~~~~~~~~~ $\widetilde{\bq}'_j = \bv_j - \widetilde{\bQ}{[:,1:j-1]} \cdot \widetilde{\bR}{[1:j-1,j]}$.} \COMMENT{macheps: $u_{crs}$}
		\STATE{4. Sketch $\widetilde{\bs}_j' = \bp_j - \widetilde{\bS}{[:,1:j-1]} \cdot \widetilde{\bR}{[1:j-1,j]}$}  \COMMENT{macheps: $u_{fine}$}
		\STATE{5. Compute the sketched norm $\widetilde{\bR}[j,j] = \|\widetilde{\bs}'_j\|$.} \COMMENT{macheps: $u_{fine}$}
		\STATE{6. Scale vector $\widetilde{\bs}_j = \widetilde{\bs}'_j/ \widetilde{\bR}[j,j] $.} \COMMENT{macheps: $u_{fine}$}
		\STATE{7. Scale vector $\widetilde{\bq}_j = \widetilde{\bq}'_j/ \widetilde{\bR}[j,j]$.} \COMMENT{macheps: $u_{fine}$}
		\ENDFOR
	\end{algorithmic}
\end{algorithm}

\begin{algorithm}[h] \caption{Steps 4-5 of RBGS: RCholeskyQR}  \label{alg:CholQR}
	\begin{algorithmic}
		\STATE{\textbf{Given:} $\bQ'_{(i)}$, $\bTheta$.} 
		\STATE{\textbf{Output}: QR fact. $\bQ'_{(i)} = \bQ_{(i)}\bR_{(i,i)}$,{where $\bR_{(i,i)}$ is upper triangular} and $\bQ_{(i)}$ is orthonormal with respect to $\langle \bTheta \cdot, \bTheta \cdot \rangle$}; a sketch $\bS_{(i)} = \bTheta \bQ_{(i)}$.
		~~\STATE{\hspace{1em}4. Compute $\bS'_{(i)}= \bTheta \bQ'_{(i)}$.}
		\STATE{\hspace{2.22em}Compute $\bR_{(i,i)}$ as the R factor of $\ell_2$-QR factorization of $\bS'_{(i)}$.}
		\STATE{\hspace{2.22em}Compute $\bQ_{(i)} = \bQ'_{(i)} (\bR_{(i,i)})^{-1}$ with  forward substitution.}
		\STATE{~~~5. Calculate the sketch of $\bQ_{(i)}$: $\bS_{(i)} = \bTheta \bQ_{(i)}$.}	
	\end{algorithmic}
\end{algorithm}

\subsubsection{Implicit RCholeskyQR factorization} \label{interQR1}
The efficiency of \Cref{alg:CholQR} on distributed architectures can be improved by using the knowledge on how $\bQ'_{(i)}$ was computed in the RBGS algorithm. This in particular can help to overcome (or to postpone) one multiplication with $\bTheta$, as described by~\cref{alg:sRBGS}. The numerical stability of such steps 4 and 5 of RBGS follows by induction from the stability of the regular RCholeskyQR, as explained in~\cref{stability}. It is then noticed that in RBGS that utilizes \cref{alg:sRBGS}  the sketch of $\bQ_{(i)}$ in step 5 of iteration $i$ can be computed together with the sketch of $\bW_{(i+1)}$ in step 1 of iteration $i+1$ with just one global synchronization between distributed processors. In this way, the resulting RBGS process would require only one global synchronization per iteration and therefore belongs to the ``one-synchronization'' family of BGS algorithms.   

\begin{algorithm}[h] \caption{Steps 4-5 of RBGS: RCholeskyQR with postponed sketching step} \label{alg:sRBGS}
	\begin{algorithmic}
		\STATE{\textbf{Given:} $\bQ'_{(i)}$, $\bTheta$, and quantities $\bS_{(1:i-1)}$, $\bR_{(1:i-1,i)}$, $\bP_{(i)}$.} 
		\STATE{\textbf{Output}: QR fact. $\bQ'_{(i)} = \bQ_{(i)}\bR_{(i,i)}$, {where $\bR_{(i,i)}$ is upper triangular} and $\bQ_{(i)}$ is orthonormal with respect to $\langle \bTheta \cdot, \bTheta \cdot \rangle$;
			a sketch $\bS_{(i)} = \bTheta \bQ_{(i)}$.}
		~~\STATE{\hspace{1em}4. Compute  $\bS'_{(i)} = \bP_{(i)} - \bS_{(1:i-1)} \bR_{(1:i-1,i)}$.}
		\STATE{\hspace{2.22em}Compute $\bR_{(i,i)}$ as the R factor of $\ell_2$-QR factorization of $\bS'_{(i)}$.}
		\STATE{\hspace{2.22em}Compute $\bQ_{(i)} = \bQ'_{(i)} (\bR_{(i,i)})^{-1}$ with  forward substitution.}
		\STATE{~~~5. Calculate the sketch of $\bQ_{(i)}$: $\bS_{(i)} = \bTheta \bQ_{(i)}$.}	
	\end{algorithmic}
\end{algorithm}

\begin{remark}[Relation between RBGS of $\bW$ and RCholeskyQR of $\bW$]
	Notice that the RBGS algorithm utilizing \cref{alg:sRBGS} in steps 4-5 computes the $\bQ$ factor as $\bW \bR^{-1}$ with forward substitution. This computation implies a strong connection between RBGS of $\bW$ and RCholeskyQR of $\bW$. The difference is that in RBGS, the R factor depends on the computed columns of the Q factor, while in RCholeskyQR it is computed solely from the sketch of $\bW$.
	
	Moreover, as noticed in~\parencite{balabanov2022randomized}, RBGS would become numerically equivalent to RCholeskyQR if in step 5 of~\cref{alg:sRBGS} the sketch of $\bQ_{(i)}$ would be obtained as $\bS_{(i)} = \bS'_{(i)} \bR_{(i,i)}^{-1}$ instead of $\bS_{(i)} = \bTheta \bQ_{(i)}$. Indeed, according to~\parencite{balabanov2022randomized}, the RBGS algorithm with such a modified step 5 can be seen as a RCholeskyQR where the $\bR_{(i,i)}$ factor is computed by BGS orthogonalization of $\bS'_{(i)}$, and the multiplication by $\bR_{(i,i)}^{-1}$ is performed block columnwise. Such RCholeskyQR and RBGS have similar computational costs in terms of flops, memory, and parallelization. However, RBGS has two advantages. First, it should be more stable in practice as it provides an R factor that accounts for errors made during the computation of $\bQ$. Second, RBGS allows for efficient certification of the solution at each iteration of the algorithm. This certification can serve, for instance, as a criterion for re-orthogonalization or a restart of a Krylov solver.
\end{remark}

\subsubsection{$\ell_2$-QR+RCholeskyQR factorization}
{Although the RCholeskyQR and RGS for inter-block orthogonalization already can provide great stability of RBGS algorithm (as shown in~\cref{stability}), this stability can be improved even more by running~\cref{alg:CholQR} twice, or by combining it with modern efficient routines for $\ell_2$-orthogonalization as described next. Note that this should have only a minor impact on the overall cost of RBGS in standard sequential architecture but not in parallel.
}

The idea here is to improve the stability of RCholeskyQR by pre-processing the matrix $\bQ'_{(i)}$ with classical routines for $\ell_2$-orthogonalization, as shown in~\cref{alg:TS-CholQR}. In principle, the pre-processing step can be done with any suitable routine for $\ell_2$-orthogonalization of tall-and-skinny matrices such as the classical or modified Gram-Schmidt with $l$ re-orthogonalizations, Householder algorithm, the tall-and-skinny QR from~\parencite{demmel2012communication} (TSQR) or any others. In each particular situation, the most suitable routine  should be selected depending on the computational architecture and programming environment. For instance, because of their popularity and reliability, the Gram-Schmidt and Householder algorithms are often available {in scientific libraries} as greatly optimized high-level routines. On the other hand, TSQR is favorable in massively parallel environments, since it reduces the amount of messages/synchronizations between processors.    

\begin{algorithm}[h] \caption{Steps 4-5 of RBGS: $\ell_2$-QR+RCholeskyQR factorization of $\bQ'_{(i)}$} \label{alg:TS-CholQR}
	\begin{algorithmic}
		\STATE{\textbf{Given:} $\bQ'_{(i)}$, $\bTheta$.} 
		\STATE{\textbf{Output}: QR fact. $\bQ'_{(i)} = \bQ_{(i)}\bR_{(i,i)}$, {where $\bR_{(i,i)}$ is upper triangular} and $\bQ_{(i)}$ is orthonormal with respect to $\langle \bTheta \cdot, \bTheta \cdot \rangle$}; a sketch $\bS_{(i)} = \bTheta \bQ_{(i)}$.
		~~\STATE{\hspace{1em}4. Compute $\ell_2$-QR fact. of $\bQ'_{(i)}$:  $\bQ'_{(i)} = \bQ^* \bR'$.}
		\STATE{\hspace{2.22em}Use~\cref{alg:CholQR}, providing $\bQ^*$ as $\bQ'_{(i)}$,  to compute RCholeskyQR $\bQ^*= \bQ_{(i)} \bR''$.}
		\STATE{\hspace{2.22em}Calculate $\bR_{(i,i)} = \bR'' \bR'$.}
		\STATE{~~~5. Calculate the sketch of $\bQ_{(i)}$: $\bS_{(i)} = \bTheta \bQ_{(i)}$.}	
	\end{algorithmic}
\end{algorithm}

It can be shown that \cref{alg:TS-CholQR} produces a stable QR factorization satisfying~\cref{eq:sRGS,eq:intercond}, if it is computed with roundoff $u = u_{fine} = F(m,n) u_{crs}$, where $F(m,n)$ is some low-degree polynomial. 

\section{Stability analysis} \label{stability}

In this section we provide a rigorous stability  analysis of the RBGS algorithm. In particular, we show that the RBGS algorithm has similar a priori as well as a posteriori stability guarantees as its single-vector counterpart (see \parencite[Section 3]{balabanov2021randomizedGS}). 

\subsection{Assumptions}
Our analysis will be based on the following assumptions. We first assume that $\bTheta$ has a bounded norm as in~\cref{eq:Thetanorm}. 
This condition is satisfied with probability at least $1-\delta$, if $\bTheta$ is $(1/2, \delta/n, 1)$ oblivious subspace embedding (see~\cref{prop:Cn2sqrtn}).  
Let the matrix $\bDelta \bQ'_{(i)}$  define the rounding error in step 3:
\begin{equation} \label{eq:DeltaQ}
\bDelta \bQ'_{(i)} = \bhQ'_{(i)} - \left (\bhW_{(i)} - \bhQ_{(1:i-1)} \bhR_{(1:i-1,i)}   \right ).  
\end{equation}
Then, the standard worst-case scenario rounding analysis gives~\cref{eq:numepsembedding2}.  Furthermore, we assume that $\bTheta$ also satisfies~\cref{eq:numepsembedding3} which, according to~\cref{thm:thetadeltax2}, holds under the {probabilistic} rounding model with probability at least $1-4\delta$, if $\bTheta$ is  $(1/8, m^{-1}\binom{n}{d}^{-1}\delta, d)$ oblivious subspace embedding, {with $d = \mathcal{O}(\log(m/\delta))$.} In its turn this condition is met by  Rademacher matrices with {$k = \mathcal{O}(\log(n)\log(m/\delta))$} rows or SRHT with $k = \mathcal{O}(\log^2(n)\log^2(m/\delta))$ rows.

Finally, it is assumed that in steps 4--5, the QR factorization and the sketch of the Q factor satisfies~\cref{eq:sRGS+,eq:intercond+} that can be {attained} by the algorithms from~\cref{interQR}. {The only exception is the implicit RCholeskyQR factorization from~\cref{interQR1}, which does not directly satisfy~\cref{eq:intercond+} and therefore requires slightly modified stability analysis of RBGS than other cases.  Fortunately the stability guarantees of the version of RBGS that uses the implicit RCholeskyQR follow directly from the guarantees of the version of RBGS that uses the regular RCholeskyQR, as is shown next. We first notice that perturbing $\bhQ'_{(i)}$ by some matrix that has (sketched) Frobenius norm $<u_{crs} F(m) \| \bW \|_\Frob$ will not change the stability guarantees from~\cref{stabgar} and their proofs up to constants. Therefore in the RBGS based on the regular RCholeskyQR we can perturb $\bhQ'_{(i)}$ to $\bhQ''_{(i)}$ so that this matrix has the sketch $\bTheta \bhQ''_{(i)}=\bhP_{(i)} - \bhS_{(1:i-1)} \bhR_{(1:i-1,i)}$. Then it is noticed that at the last line of step 4 of the RCholeskyQR we can perturb the matrix $\bhQ''_{(i)}$  back to $\bhQ'_{(i)}$ as this can increase $\mathrm{cond}(\bTheta \bhQ_{(i)})$ only by $\mathcal{O}(u_{crs} F(m)  \| \bW \|_\Frob \|\bR_{(i,i)}\|^{-1}_\Frob )= \mathcal{O}(u_{crs} \mathrm{cond}(\bW) F(m))$, which can be shown straightforwardly by following the proof of stability of RCholeskyQR in~\parencite{balabanov2022randomized}. The argument is finished by noticing that the implicit RCholeskyQR can be viewed exactly as the regular RCholeskyQR with such two perturbations.}

\begin{assumptions} \label{thm:Thetaasmpts}
	It is assumed that for some $\varepsilon \leq 1/2$, 
	\begin{equation} \label{eq:Thetanorm}
	\| \bTheta \|_\Frob \leq \sqrt{1+\varepsilon} \sqrt{n}. 
	\end{equation}
	Furthermore,  we assume that
	\begin{subequations} 
		\begin{align}
		|\bDelta {{\bQ}'_{(i)}}|	&\leq  1.02 u_{crs} (| \bhW_{(i)} | +  i m_p | \bhQ_{(1:i-1)}| |\bhR_{(1:i-1,i)} |) \label{eq:numepsembedding2} \\
		\|\bTheta \bDelta {{\bQ}'_{(i)}} \|_\Frob &\leq 1.02 u_{crs} \sqrt{1+\varepsilon} \|| \bhW_{(i)} | +  i m_p | \bhQ_{(1:i-1)}| |\bhR_{(1:i-1,i)} |\|_\Frob \label{eq:numepsembedding3} 
		\end{align}
	\end{subequations}
	with $1 \leq i \leq p$ and $1 \leq t \leq m_p$.
	Finally, it is assumed that in steps 4--5 of~\cref{alg:RBGS}, we have
	\begin{subequations} \label{eq:sRGS+}
		\begin{align}
		&|{\bhQ}'_{(i)} - {\bhQ}_{(i)} \bhR_{(i,i)}|  \leq u_{fine} m (|{\bhQ}'_{(i)}|+|{\bhQ}_{(i)}| | \bhR_{(i,i)} |) \leq 0.1 u_{crs} |{\bhQ}_{(i)}| | \bhR_{(i,i)} | \label{eq:sRGS1+} \\
		\begin{split}
		&\|\bTheta({\bhQ}'_{(i)} - {\bhQ}_{(i)} \bhR_{(i,i)}) \|_\Frob  \leq u_{fine} m \|\bTheta \|_\Frob \||{\bhQ}'_{(i)}|+|{\bhQ}_{(i)}| | \bhR_{(i,i)} |\| \\
		& ~~~~~~~~~~~~~~~~~~~~~~~~~~~~~\leq  0.1 u_{crs} \| {\bhQ}_{(i)} \|  \| \bhR_{(i,i)}\|,
		\end{split}\label{eq:sRGS2+}
		\end{align}
	\end{subequations}
	and, if $\bTheta$ is {a} $\varepsilon'$-embedding for $\bhQ'_{(i)}$, 
	\begin{equation} \label{eq:intercond+}
	1-0.1 u_{crs} \cond(\bhQ'_{(i)}) \leq  \sigma_{min}(\bTheta \bhQ_{(i)}) \leq \sigma_{max}(\bTheta \bhQ_{(i)}) \leq 1+0.1 u_{crs} \cond(\bhQ'_{(i)}), 
	\end{equation}
	and $\| \bhQ_{(i)}\|\leq (1-\varepsilon')^{-1/2} \|\bTheta \bhQ_{(i)}\|$.
\end{assumptions}	

\subsection{Stability guarantees of RBGS algorithm} \label{stabgar}

Our stability  analysis will rely on the condition that $\bTheta$ satisfies the $\varepsilon$-embedding property for $\bhW$ and $\bhQ$. See~\cref{epsilonprop} for a  characterization of this property.

\subsubsection{A posteriori analysis of RBGS algorithm}
Let us first give an a posteriori characterization of the RBGS algorithm. Such characterization can be performed by measuring or bounding coefficients $\Delta^{(p)}  = \|\bI - \bhS^\mathrm{T} \bhS \|_\Frob$ and $\tilde{\Delta}^{(p)}= \frac{\|\bhP - \bhS \bhR \|_\Frob}{\|\bhP\|_\Frob}$, as done in~\parencite{balabanov2021randomizedGS}. We have the following result, which {is an analogue of~\parencite[Theorem 3.2]{balabanov2021randomizedGS} for RBGS}.
\begin{theorem}\label{thm:maintheorem1}
	Consider~\Cref{alg:RBGS}. {Assume that} 
	$$ {100} m^{1/2} n^{3/2} u_{fine} \leq u_{crs} \leq 0.01,$$
	{along with~\cref{thm:Thetaasmpts}, possibly excluding~\cref{eq:intercond+}.}	
	
	{If $\bTheta$ is {an} $\varepsilon$-embedding for $\bhQ$ and $\bhW$, with $\varepsilon \leq 1/2$, and if $\Delta^{(p)},\tilde{\Delta}^{(p)} \leq 0.1$,}
	then the following inequalities hold:			
	\begin{subequations}
		\begin{equation*}
		{
			(1+\varepsilon)^{-1/2}(1-{\Delta^{(p)}}- 0.1 u_{crs}) \leq \sigma_{min}(\bhQ)   \leq \sigma_{max}(\bhQ) \leq (1-\varepsilon)^{-1/2}(1+{\Delta^{(p)}}+ 0.1 u_{crs})}
		\end{equation*}
		\normalsize
		\begin{equation*}
		\|\bhW - \bhQ \bhR\|_\Frob \leq  4 u_{crs} m^{3/2} \|\bhW\|_\Frob. 
		\end{equation*}
	\end{subequations}
	
	\begin{proof}
		See~\cref{proofs}.
	\end{proof}
	
\end{theorem}
\begin{remark} \label{thm:maintheorem1+}
	In~\cref{thm:maintheorem1}, we also have
	$$\|\bTheta (\bhW - \bhQ \bhR)\|_\Frob \leq  5 u_{crs} m^{3/2} \|\bhW\|_\Frob.$$
\end{remark}

\Cref{thm:maintheorem1} implies the numerical stability of the RBGS algorithm if $\Delta^{(p)}$ and $\tilde{\Delta}^{(p)}$ are $\leq 0.1$. These coefficients can be efficiently computed a posteriori from the sketches $\bhS$ and $\bhP$ and the R factor $\bhR$, thus providing a way for the certification of the solution. Such certification, in particular, does not involve any assumptions on $\cond(\bhW)$, the accuracy of the least-squares solution in step 2, and the stability of inter-block orthogonalization in steps 4--5. Furthermore, we would like to highlight the very mild condition $u_{crs} = \mathcal{O}(m^{-3/2})$ on the working (coarse) unit roundoff to guarantee {the} accuracy of the algorithm, which is in particular independent of the high-dimension $n$.

\subsubsection{A priori analysis of RBGS algorithm}
Clearly, to get a priori bounds for {$\Delta^{(p)}  = \|\bI - \bhS^\mathrm{T} \bhS \|_\Frob$ and $\tilde{\Delta}^{(p)}= \frac{\|\bhP - \bhS \bhR \|_\Frob}{\|\bhP\|_\Frob}$} we need more assumptions than {those stated in}~\cref{thm:maintheorem1}.  In particular, it is necessary to impose a stability condition on the least squares solver used in step 2. We also need $\bhW$ to be numerically {full rank}, i.e., to satisfy $u_{crs} \leq \mathcal{O}(\cond(\bhW)^{-1})$.

The following a priori guarantee of stability of {the} RBGS algorithm is, in fact, an analogue of~\parencite[Theorem 3.3]{balabanov2021randomizedGS} for {the} single-vector RGS algorithm.

\begin{theorem} \label{thm:maintheorem2}
	Consider~\cref{alg:RBGS} with a backward-stable solver (e.g., based on Richardson iterations) satisfying~\cref{thm:asssolver}.
	
	{Under~\cref{thm:Thetaasmpts},} assume that $\bTheta$ is {an} $\varepsilon$-embedding for $\bhQ_{(1:p-1)}$  and $\bhW$, with $\varepsilon \leq 1/2$.  If	$$u_{crs} \leq 10^{-3} \cond(\bhW)^{-1} m^{-2} \textit{ and } u_{fine} \leq 10^{-2} m^{-1/2} n^{-3/2} u_{crs},$$
	then $\Delta^{(p)}$ and $\tilde{\Delta}^{(p)}$ are bounded by
	\begin{align} 
	\tilde{\Delta}^{(p)}
	& \leq 4.2 u_{crs}  m^{3/2} \|\bhW\|_\Frob/\|\bhP\|_\Frob \leq 6 u_{crs} m^{3/2}, \label{eq:maintheorem21} \\
	\Delta^{(p)}
	&\leq 20 u_{crs}   m^{2} \cond{(\bhW)}. \label{eq:maintheorem22}
	\end{align}
	\begin{proof}
		See~\cref{proofs}.						
	\end{proof}
\end{theorem}
\Cref{thm:maintheorem2} states that the RBGS algorithm is stable unless the input matrix $\bhW$ is numerically {rank-deficient}. 
{This stability guarantee is seen in other stable deterministic algorithms such as MGS, CGS2, BCGS2, and others.}
Furthermore, the stability is proven for {the} {working} unit roundoff independent of {the high dimension $n$.} {This unique feature of randomized algorithms can be especially interesting for large-scale problems solved on low-precision arithmetic architectures.}

\subsection{Epsilon embedding property} \label{epsilonprop}
The stability analysis in~\cref{stability} holds if $\bTheta$ satisfies the $\varepsilon$-embedding property for $\bhQ$ and $\bhW$. In this section we analyze this property.

We consider the case when $\bhW$ and $\bTheta$ are independent of each other. Then, if $\bTheta$ is {a} $(\varepsilon, \delta, m)$ oblivious $\ell_2$-subspace embedding, it satisfies the $\varepsilon$-embedding property for $\bhW$ with high probability. Below, we show that{,} in this case $\bTheta$ will also satisfy {the} $\varepsilon$-embedding property for $\bhQ$ with moderately increased value of $\varepsilon$. This result is basically the RBGS counterpart of~\parencite[Proposition 3.6]{balabanov2021randomizedGS}. 
\begin{proposition}\label{thm:aprioribound}
	Consider~\cref{alg:RBGS} with a backward-stable solver satisfying~\cref{thm:asssolver},
	$$u_{crs} \leq 10^{-3} \cond(\bhW)^{-1} m^{-2} \textit{ and } u_{fine} \leq 10^{-2} m^{-1/2} n^{-3/2} u_{crs}.$$
	{Under~\cref{thm:Thetaasmpts},}  if $\bTheta$ is {a} $\varepsilon$-embedding for $\bhW$, with $\varepsilon \leq 1/4$, then it satisfies the $\varepsilon'$-embedding property for $\bhQ$ with 
	$\varepsilon' = 2 \varepsilon + 180 u_{crs} m^2 \cond(\bhW). $
	\begin{proof}  
		See~\cref{proofs}.						
	\end{proof}
\end{proposition}

The $\varepsilon$-embedding property will likely hold even when matrix $\bhW$ depends on $\bTheta$. In such a case the quality of $\bTheta$ can be certified a posteriori by computing additional sketches $\bPhi \bhQ$ and $\bPhi \bhW$, associated with a new sketching matrix $\bPhi$, in addition to the sketches $\bTheta \bhQ$ and $\bTheta \bhW$.  Then one may characterize the quality of $\bTheta$ by measuring the orthogonality of $(\bTheta \bhW) \bX$ and $(\bTheta \bhQ)  \bY$, where $\bX$ and $\bY$ are inverses of R factors of $\bPhi \bhW$ and $\bPhi \bhQ$, as is described in~\parencite[Propositions 3.6-3.8]{balabanov2021randomizedGS} extrapolated from~\parencite{balabanov2019randomized2}. For efficiency in terms of cache or communication, at each iteration, the products $\bPhi \bhW_{(i)}$ and $\bPhi \bhQ_{(i)}$  can be computed together with $\bTheta \bhW_{(i)}$ in step 2, and  $\bTheta \bhQ_{(i)}$ in steps 4-5, respectively.

\section{Randomized block Krylov methods} \label{applications}

In this section we discuss practical applications of the methodology. We particularly focus on {improving the} efficiency of popular block Krylov methods, such as the block GMRES and FOM, and the RR method, for solving block linear systems of the form $\bA \bX = \bB$
and eigenvalue problems of the form $\bA \bX = \bX \bLambda$,
where $\bA$ is large, possibly non-symmetric $n \times n$ matrix, $\bB$  is $n \times m_p$ matrix, and $\bLambda$ is $m_p \times m_p$ diagonal matrix of the extreme eigenvalues of $\bA$. 

The block Krylov  methods proceed with approximation of $\bX$ by a projection $\bX^{(j)}$ onto {the} Krylov space $\mathcal{K}^{(j)}(\bA, \bB)$, {which is} defined as
$$\mathcal{K}^{(j)}(\bA, \bB) := \mathrm{span} \{\bB, \bA \bB, \hdots, \bA^{j-1} \bB \},$$
with $j$ being the order of the subspace.

The GMRES method computes $\bX^{(j)}$ that minimizes the Frobenius norm of the residual, while the FOM and RR methods seek $\bX^{(j)}$ that is optimal in the Galerkin sense. Both kinds of methods first construct an orthonormal basis for $\mathcal{K}^{(p)}(\bA, \bB)$ with  GS orthogonalization, called Arnoldi iteration, and then determine the coordinates of the columns of $\bX^{(j)}$ in the computed basis.

\subsection{Krylov basis computation: randomized block Arnoldi iteration} \label{RBAi}
The Arnoldi algorithm produces orthonormal matrix $\bQ = \bQ_{(1:p)}$ satisfying the Arnoldi identity 
$ \bA \bQ_{(1:p-1)} =  \bQ_{(1:p)} \bH, $
where $\bH = \bH_{(1:p,1:p-1)}$ is block upper {Hessenberg} matrix. The Arnoldi algorithm can be viewed as a block-wise QR factorization of matrix $[\bB, \bA  \bQ_{(1:p-1)}]$. In this context, the {Arnoldi} matrix $\bH$ can be seen as the R factor $\bR = \bR_{(1:p,1:p)}$ without the first column of subblocks, that is, as $\bR_{(1:p,2:p)}$.

Below, we propose a randomized Arnoldi process based on {the} RBGS algorithm. Note that unlike standard methods,~\cref{alg:RGS-Arnoldi} produces {a} Krylov basis orthonormal with respect to the sketched product $\langle \bTheta \cdot, \bTheta \cdot \rangle$. 

\begin{algorithm}[h] \caption{RBGS-Arnoldi algorithm} \label{alg:RGS-Arnoldi}
	\begin{algorithmic}
		\STATE{\textbf{Given:} $n \times n$ matrix  $\bA$, $n \times m_p$ matrix $\bB$, $k \times n$ matrix $\bTheta$ with $k \ll n$, parameter $p$.}
		\STATE{\textbf{Output}: $n \times m$ factor $\bQ = \bQ_{(1:p)}$ and $m \times m$ upper triangular factor $\bR = \bR_{(1:p,1:p)}$.}
		\STATE{1. Set  $\bW_{(1)} = \bB$ and perform $1$-st iteration of {the} RBGS algorithm.}
		\FOR{ $i=2:p$} 		
		\STATE{2. Compute $\bW_{(i)} = \bA \bQ_{(i-1)}$. }
		\STATE{3. Perform $i$-th iteration of {the} RBGS algorithm.}
		\ENDFOR
		\STATE{4. (Optional) compute $\Delta^{(p)}$ and $\tilde{\Delta}^{(p)}$.  Use~\cref{thm:maintheorem1} to certify the output.}
	\end{algorithmic}
\end{algorithm}

In step~2 of~\cref{alg:RGS-Arnoldi}, the computation of the matrix-vector product can be executed either with roundoff $u_{fine}$ or $u_{crs}$ depending on the situation. 
The stability guarantees of~\cref{alg:RGS-Arnoldi} can be obtained directly from~\cref{thm:maintheorem1,thm:maintheorem2} with standard stability analysis similar to that from~\parencite[Section 4.1]{balabanov2021randomizedGS}. In particular, it can be shown that, under the stability conditions of RBGS, the computed $\bhQ$ and $\bhH$ satisfy 
\begin{equation} \label{eq:pArnoldi}
(\bA+\bDelta\bA) \bhQ_{(1:p-1)} = \bhQ_{(1:p)} \bhH,
\end{equation}
with $\|\bDelta\bA\|$ close to machine precision, and $\cond(\bhQ) = \mathcal{O}(1)$. We leave the precise analysis of this fact outside of the scope of this manuscript.

\subsection{Linear systems: randomized block GMRES method} \label{RBGSRES}

Let us now discuss {solving} block linear systems $\bA \bX = \bB$ with GMRES.

The GMRES method computes the approximate solution $\bX^{(p-1)} = \bU = [\bu_1, \hdots, \bu_{m_p}]$:
\begin{equation} \label{eq:GMRESsol}
\bu_i = \bhQ_{(1:p-1)} \arg \min_\bz \| \bhH \bz - \bhR_{(1:p,1)} \be_i \|, 
\end{equation}
performed in sufficient precision, where $\be_i$ is the $i$-th column of the identity matrix. 
Let us now characterize quasi-optimality of such a projection when $\bhQ$ and $\bhH$ were obtained with the RBGS-Arnoldi algorithm. 
Under the stability conditions of RBGS,  the Arnoldi identity~\cref{eq:pArnoldi} implies that
$$ \| (\bA+\bDelta\bA)\bhQ_{(1:p-1)} \bz - \bb_i \| =  \| \bhQ_{(1:p)}(\bhH \bz- \bhR_{(1:p,1)}\be_i)  \| \leq \| \bhQ \| \| \bhH \bz - \bhR_{(1:p,1)} \be_i \|,$$
and, similarly,
$$ \| (\bA+\bDelta\bA) \bhQ_{(1:p-1)} \bz - \bb_i \| \geq \sigma_{min} (\bhQ) \| \bhH \bz - \bhR_{(1:p,1)} \be_i \|.$$
These two relations imply that
$$ \|(\bA+\bDelta\bA) \bu_i - \bb_i \| \leq \cond(\bhQ) \min_{\bv \in Q_{p-1}} \| (\bA+\bDelta\bA) \bv- \bb_i \|,  $$
with $Q_{p-1} = \mathrm{range}(\bhQ_{(1:p-1)}) = \mathcal{K}_{p-1}(\bA+\bDelta\bA,\bB).$ Consequently, the randomized version of GMRES provides a solution which minimizes the norm of the residual associated with a slightly perturbed matrix, up to a factor of order 1.

\subsection{Linear systems: randomized block FOM method} \label{rmk:RFOM}
In contrast to GMRES, the {FOM method} obtains solution $\bX^{(p-1)}=\bU=[\bu_1, \hdots, \bu_{m_p}]$ by imposing a Galerkin orthogonality condition on the residuals (in exact arithmetic):
\begin{equation}\label{eq:GalerkinFOM}
	\langle \bv, \br(\bu_i,\bb_i) \rangle = 0,~~~ \forall \bv \in Q_{p-1},~1\leq i \leq m_p, 
\end{equation}
where $Q_{p-1} = \mathrm{range}(\bQ_{(1:p-1)})$ and $\br(\bu_i,\bb_i) = \bA \bu_i  - \bb_i$ is the residual associated with the $i$-th right-hand-side. The Galerkin projection may be a more appropriate choice than {the minimal-residual projection provided by GMRES}, when the quality of the solution is measured with an  energy error rather than the residual error.  When $\bQ_{(1:p-1)}$ is obtained with the {traditional} Arnoldi method, the solution to~\cref{eq:GalerkinFOM} is given by 
\begin{equation}\label{eq:RFOMsol}
	\bU = \bQ_{(1:p-1)} (\bH_{(1 : p-1,1:p-1)})^{-1} \bR_{(1:p-1,1)}.
\end{equation}

Suppose now that $\bQ_{(1:p-1)}$ was obtained by the randomized Arnoldi method. Then, by noticing that 
$$\bH_{(1 : p-1,1:p-1)} = [\bI~\bnull_{(1:p-1,p)}] \bH = (\bTheta\bQ_{(1:p-1)})^\mathrm{T} \bTheta \bQ_{(1:p-1)} \bH = (\bTheta \bQ_{(1:p-1)})^\mathrm{T} \bTheta \bA \bQ_{(1:p-1)}, $$ 
where $\bnull_{(1:p-1,p)}$ is a null matrix of size of $\bH_{(p,1:p-1)}^\mathrm{T}$,  we deduce that the solution~\cref{eq:RFOMsol} satisfies 
\begin{equation}\label{eq:skGalerkinFOM}
	\langle \bTheta \bv, \bTheta \br(\bu_i,\bb_i) \rangle = 0,~~~ \forall \bv \in Q_{p-1},~1\leq i \leq m_p. 
\end{equation} 
The relation~\cref{eq:skGalerkinFOM} can be viewed as the sketched version of the Galerkin orthogonality condition.  To our knowledge, it was first used in model order reduction community  to obtain a reduced-basis solution of parametric linear systems~\parencite{balabanov2019randomized}. By using similar considerations as in~\parencite{balabanov2019randomized}, one can show that~\cref{eq:skGalerkinFOM} preserves the quality of the classical Galerkin projection when $\bTheta$ is a $\varepsilon$-embedding for $Q_{p}$ with $\varepsilon\cond(\bA)<1$, though, in practice, this condition can be too pessimistic~\parencite{balabanov2019randomized}.  We leave the further development of the randomized FOM method for future research.

\begin{remark}
	{It should be noted that the reasoning of~\cref{RBGSRES,rmk:RFOM} could be reversed. We could first derive the sketched {minimal-residual} projection
		$\bu_i = \arg \min_{\bv \in Q_{p-1}} \allowbreak \|\bTheta \br(\bv,\bb_i) \|$ and sketched Galerkin orthogonality condition~\cref{eq:skGalerkinFOM} by 
		replacing the $\ell_2$-inner products and $\ell_2$-norms in the standard {minimal-residual} equation $\bu_i = \arg \min_{\bv \in Q_{p-1}} \allowbreak\|\br(\bv,\bb_i) \|$  and the standard Galerkin orthogonality condition~\cref{eq:GalerkinFOM} by sketched ones, similarly as was done in~\parencite{balabanov2019randomized,balabanov2019randomized2}. After that, it could be realized that the solution to the sketched {minimal-residual} and Galerkin equations can be obtained by orthogonalizing the Krylov basis  with respect to  $\langle \bTheta \cdot , \bTheta \cdot \rangle$ and using the classical identities~\cref{eq:GMRESsol,eq:RFOMsol}.}
\end{remark}

\subsection{Eigenvalue problems: randomized RR method} \label{rRR}

Next, we consider the computation of the extreme eigenpairs of $\bA$. Note that the methodology proposed below can be readily used for finding eigenpairs in the desired region by introducing a shift to the matrix. In
addition, our methodology can be extended to inverse methods by replacing $\bA$ with $\bA^{-1}$.
In this case, the application of the inverse  can be readily  performed with (possibly preconditioned) randomized block GMRES method from~\cref{RBGSRES} or other efficient methods.

To simplify the presentation, the analysis of the RR method and its randomized version will be  provided only in infinite precision arithmetic. The stability of these approaches follows directly from the stability of the Arnoldi iteration.

Let $\bQ = \bQ_{(1:p)}$ be {a} Krylov basis generated with the standard or randomized block-Arnoldi algorithm.  The standard RR method, based on the Arnoldi iteration,   approximates the extreme eigenpairs $(\lambda, \bx)$ of $\bA$  by 
\begin{equation} \label{eq:gensol}
(\lambda, \bx) \approx (\mu, \bu) = (\mu, \bQ_{(1:p-1)} \by),
\end{equation}
where $(\mu, \by)$ are the corresponding extreme eigenpairs of $\bH_{(1 : p-1,1:p-1)} $. Furthermore the residual error of an approximate eigenpair  $(\mu, \bu)$ is estimated by $ \| \bH_{(p,1:p-1)} \by \|$.

For the standard block-Arnoldi algorithm, we have 
{$$\bH_{(1 : p-1,1:p-1)} = [\bI~ \bnull_{(1:p-1,p)}] \bH = (\bQ_{(1: p-1)})^\mathrm{T} \bQ \bH = (\bQ_{(1: p-1)})^\mathrm{T} \bA \bQ_{(1: p-1)},$$ }
and therefore 
\begin{equation} \label{eq:Galerkin}
\langle \bv, \br(\bu,\mu \bu) \rangle  = 0,~~~ \forall \bv \in Q_{p-1}, 
\end{equation}
where $Q_{p-1} := \mathrm{range}(\bQ_{(1: p-1)})$ and $\br(\bu,\mu \bu) = \bA \bu  - \mu \bu$ is the residual associated with $(\mu, \bu)$. The relation~\cref{eq:Galerkin} is known as the Galerkin orthogonality condition. At the same time, we have
\begin{subequations} \label{eq:RRresbound}
	\begin{align}
	\|\br(\mu,\bu) \| \geq  \sigma_{min}(\bQ)
	\left \|\bH \by - \mu \left [\begin{matrix} \by \\ \bnull \end{matrix} \right ] \right \| &= \sigma_{min}(\bQ)\| \bH_{(p,1:p-1)} \by \|, \label{eq:RRresbound1} \\
	\intertext{and, similarly,} 
	\|\br(\mu,\bu) \| &\leq \sigma_{max}(\bQ) \| \bH_{(p,1:p-1)} \by \|. \label{eq:RRresbound2}
	\end{align}
\end{subequations}
Since, the classical Arnoldi iteration produces an $\ell_2$-orthogonal Q factor, hence in this case, the quantity $\| \bH_{(p,1:p-1)} \by \|$ represents exactly the residual error of $(\mu, \bu)$.

On the other hand, the randomized block-Arnoldi algorithm produces a matrix $\bQ$ orthogonal with respect to $\langle \bTheta \cdot , \bTheta \cdot \rangle$. This implies that
$$\bH_{(1 : p-1,1:p-1)} = [\bI~ \bnull_{(1:p-1,p)}] \bH = (\bTheta\bQ_{(1: p-1)})^\mathrm{T} (\bTheta \bQ_{(1: p-1)}) \bH = (\bTheta \bQ_{(1: p-1)})^\mathrm{T} \bTheta \bA \bQ_{(1: p-1)}, $$
or equivalently, that $(\mu, \bu)$ satisfies the following sketched version of the Galerkin orthogonality condition:
\begin{equation} \label{eq:skGalerkin}
\langle \bTheta \bv, \bTheta \br(\bu,\mu\bu) \rangle  = 0,~~~ \forall \bv \in Q_{p-1},
\end{equation}
similar to the sketched Galerkin condition for linear systems in~\cref{rmk:RFOM}. 
Unlike for the sketched {minimal-residual} projection in GMRES, the optimality for the sketched Galerkin projection unfortunately does not trivially follow from the $\varepsilon$-embedding property of $\bTheta$. A characterization of the accuracy of this projection should be derived by reformulating the methodology in terms of projection operators similarly to state-of-the-art analysis provided, for instance, in \parencite[Section 4.3]{saad2011numerical}. In particular, we notice that the Galerkin orthogonality condition~\cref{eq:Galerkin} can be expressed as 
$$ \bPi_{Q_{p-1}} \bA \bPi_{Q_{p-1}}  \bu = \mu\bu,$$
where $\bPi_{Q_{p-1}}$ is the $\ell_2$-orthogonal projector onto $Q_{p-1}$.
In other words, the classical RR method can be interpreted as approximation of eigenpairs of $\bA$ by the eigenpairs of \emph{approximate} operator $\bPi_{Q_{p-1}} \bA \bPi_{Q_{p-1}}$. 
Similarly, the sketched Galerkin orthogonality condition~\cref{eq:skGalerkin} can be expressed as 
$$ \bPi^\bTheta_{Q_{p-1}} \bA \bPi^\bTheta_{Q_{p-1}}  \bu = \mu\bu,$$
where $\bPi^\bTheta_{Q_{p-1}} = \bQ_{(1:p-1)} (\bTheta \bQ_{(1:p-1)})^\dagger \bTheta$ is an orthogonal projector onto $Q_{p-1}$ with respect to the sketched inner product $\langle \bTheta \cdot, \bTheta \cdot \rangle$.
We see that the randomized RR method corresponds to taking the {approximate} operator as $\bPi^\bTheta_{Q_{p-1}} \bA \bPi^\bTheta_{Q_{p-1}}$ instead of classical $\bPi_{Q_{p-1}} \bA \bPi_{Q_{p-1}}$. This connection can be used to extrapolate the results from~\parencite[Section 4.3]{saad2011numerical}, such as Theorem 4.3 that bounds the residual error of the exact eigenpair with respect to the approximate operator, from classical methods to their sketched variants. For a more detailed discussion, please refer to the supplementary materials.

Note that, independently of this work, the sketched Galerkin orthogonality condition was also used by Nakatsukasa and Tropp in their recent paper~\parencite{nakatsukasa2021fast}.

Since the $\varepsilon$-embedding property of $\bTheta$ implies  that $\cond(\bQ)=1+\mathcal{O}(\varepsilon)$, hence according to~\cref{eq:RRresbound}, the quantity  $\| \bH_{(p,1:p-1)} \by \|$ estimates well the residual error.

The RR algorithm based on RBGS-Arnoldi with restarting is depicted in~\cref{alg:RGS-Arnoldi_eig_max}.
\begin{remark} \label{thm:skGalerkintoGalerkin}
	When the classical Galerkin projection is preferable to its sketched variant, say, because of its strong optimality properties for symmetric (or Hermitian) operators, the RBGS-Arnoldi algorithm can be used to obtain this projection {instead of the sketched one}.
	For this the $\ell_2$-orthogonal Krylov basis matrix $\bQ$ and the {Arnoldi} matrix $\bH$ in~\cref{eq:gensol} {can be obtained} with the RBGS-Arnoldi algorithm followed by an additional Cholesky QR step, as is depicted in~\cref{postproc}.
	In particular, this situation can be accounted for in~\cref{alg:RGS-Arnoldi_eig_max} by adding the following line between step 1 and step 2:
	{``Compute $\bQ^\mathrm{T}\bQ$ and obtain its Cholesky factor $\bR'$. Set $\bQ \leftarrow \bQ \bR'^{-1}$, $\bR \leftarrow \bR' \bR$.''}
	The cost of this additional step is dominated by the computation of product $\bQ^\mathrm{T} \bQ$ (since we can omit computing $\bQ \bR'^{-1}$), which requires similar number of flops as the entire RBGS-Arnoldi algorithm. Nevertheless, as a single BLAS3 operation, it is very well suited to modern computational architectures. 
\end{remark}

\begin{algorithm}[h] \caption{\small Rand. RR algorithm for extreme eigenpairs with restarting} \label{alg:RGS-Arnoldi_eig_max}
	\begin{algorithmic}
		\STATE{\textbf{Given:} $n \times n$, $n \times m_p$, $k \times n$ matrices  $\bA$, $\bB$, and $\bTheta$ with $m \leq  k \ll n$, param. $p$ and $N_{iter}$.}
		\STATE{\textbf{Output}: $\bLambda $ and $\bX$. }
		\FOR{$i = 1:N_{iter}$}
		\STATE{1. Perform RBGS-Arnoldi \Cref{alg:RGS-Arnoldi} returning $\bQ$ and $\bR$.  Set $\bH = \bR_{(1:p,2:p)}$.}
		\STATE{2. Compute diagonal matrix $\bLambda$ with  $m_p$ extreme eigenvalues of $\bH_{(1:p-1,1:p-1)} $ \\ on the diagonal, and the matrix  of associated eigenvectors $\bY$.}
		\STATE{3. (Optional) Compute $\| \bH_{(p,1:p-1)} \by \|$ to characterize the approximation error.}
		\STATE{4. Compute $ \bB = \bQ_{(1:p-1)}\bY$.}
		\ENDFOR
		\STATE{5. Normalize $\bB$ with respect to $\ell_2$-norm and return it as $\bX$.}
	\end{algorithmic}
\end{algorithm}

\subsection{Further applications}

\subsubsection{$s$-step Krylov methods}

The proposed RBGS algorithm can be used to improve the $s$-step Krylov methods~\parencite{hoemmen2010communication}.  For simplicity, consider the case of a linear system with only one right-hand side or an eigenvalue problem approximated by a Krylov space associated with only one generating vector, i.e., taking $\bB = \bb$. Then the goal is to compute a basis for $\mathcal{K}^{(m)}(\bA, \bb)$ satisfying the {Arnoldi} identity. With this basis, an approximate solution of a linear system or an eigenvalue problem can be obtained using the identities from~\cref{RBGSRES,rmk:RFOM,rRR}, setting $\bB = \bb$.

The $s$-step Krylov methods use efficient matrix power kernels that can output a matrix $F_{{s+1}}(\bv)$ of basis vectors for $\mathcal{K}^{(s+1)}(\bA, \bv)$, of the form
$$F_{{s+1}}(\bv) = [p_0(\bA)\bv,p_1(\bA)\bv,\hdots ,p_s(\bA)\bv], $$
where $\bv$ is some input vector, $p_0(\bA)$ is usually taken as $\bI$, and  $p_1(\bA), \hdots p_s(\bA)$ are some suitable polynomials {that aim to make $F_{{s+1}}(\bv)$ not too badly conditioned. }
Popular options for $p_0(\bA), \hdots p_s(\bA)$ are the monomials and Newton, or Chebyshev polynomials.  A relatively large $s$ {(say $s = 10$ or $30$~\parencite{hoemmen2010communication})} may lead to stability problems even when a polynomial basis is used. {Consequently, to obtain a basis for $\mathcal{K}^{(m)}(\bA, \bb)$ of dimension suitable for practical applications, it becomes necessary}  to proceed with a block-wise generation, at each iteration, computing a small block of a Krylov basis by using the matrix power kernel and subsequently orthogonalizing it against the previously computed vectors with a BGS approach.
In more concrete terms, {we have to} perform a block-wise QR factorization of the matrix $\bW = \bW_{(1:p)}$ generated recursively as 
$$ \bW_{(i)} = \left \{ \begin{matrix}   \bb &\text{ if $i=1$}\\
F_{{s+1}}^*(\bq_{s(i-2)+1}) &\text{ if $i= 2, \hdots, p$} \end{matrix} \right.$$ 
where $F_{{s+1}}^*(\bv)$ corresponds to the matrix $F_{{s+1}}(\bv)$ without the first column, and $\bq_{s(i-2)+1}$ is the last Krylov basis vector computed at iteration~$i-1$. Such block-wise orthogonalization  can readily be done with the RBGS algorithm, as is depicted in~\cref{alg:RGS-sstep-Arnoldi}, with all the computational benefits of this algorithm.

\begin{algorithm}[h] \caption{$s$-step RBGS-Arnoldi algorithm} \label{alg:RGS-sstep-Arnoldi}
	\begin{algorithmic}
		\STATE{\textbf{Given:} $n \times n$ matrix  $\bA$, $n \times 1$ vector $\bb$, $k \times n$ matrix $\bTheta$ with $k \ll n$, parameters $p$ and $s$.}
		\STATE{\textbf{Output}: $n \times m$ Krylov basis matrix $\bQ = \bQ_{(1:p)}$ and {Arnoldi} matrix $\bH$.}
		\STATE{1. Set  $\bW_{(1)} = \bb$ and perform $1$-st iteration of {the} RBGS algorithm.}
		\FOR{ $i=2:p$} 		
		\STATE{2. Compute $\bW_{(i)} = F_{{s+1}}^*(\bq_{s(i-2)+1})$ and perform $i$-th iteration of {the} RBGS algorithm. }
		\ENDFOR
		\State{3. Obtain $\bH = \bS^\dagger {(\bTheta \bA \bQ_{(1:p-1)})}$ by QR or SVD of $\bS$, or solving a small least-squares problem. (Optional) set the below-subdiagonal elements of $\bH$  to zeros.}
		\STATE{4. (Optional) compute $\Delta^{(p)}$ and $\tilde{\Delta}^{(p)}$.  Use~\cref{thm:maintheorem1} to certify the output.}
	\end{algorithmic}
\end{algorithm}

To obtain the Arnoldi matrix $\mathbf{H}$ in step 3, \Cref{alg:RGS-sstep-Arnoldi} utilizes the sketch $\mathbf{C} := \bTheta \mathbf{A} \mathbf{Q}_{(1:p-1)}$. Computing $\mathbf{C}$ should incur only minor costs in terms of flops and communication. It can be beneficial to compute $\bTheta (\mathbf{A} \mathbf{Q}_{(i-1)})$ together with $\bTheta \mathbf{W}_{(i)}$ at the $i$-th iteration of the RBGS algorithm. Alternatively, $\bC$ can be obtained through a recursive procedure without any high-dimensional operations.  At iteration $i \geq 2$, we can first compute $\mathbf{c}_{s(i-2)+1}$ and $\bTheta \mathbf{A} \overline{\mathbf{W}_{(i)}}$, where $\overline{\mathbf{W}_{(i)}}$ represents $\mathbf{W}_{(i)}$ without the last column, by using the following relation: $[\mathbf{c}_{s(i-2)+1}~ \bTheta \mathbf{A} \overline{\mathbf{W}_{(i)}}] = [\mathbf{s}_{s(i-2)+1}~ \mathbf{P}_{(i)}] \mathbf{T}_{s+1}$. Here, $\mathbf{P}_{(i)}$ and $\mathbf{s}_{s(i-2)+1}$ denote the sketches of $\mathbf{W}_{(i)}$ and $\mathbf{q}_{s(i-2)+1}$, respectively, computed during the $i$-th iteration of the RBGS algorithm, and $\mathbf{T}_{s+1}$ is the Arnoldi factor associated with the matrix power kernel $F_{{s+1}}(\mathbf{v})$, which satisfies $\mathbf{A} F_{{s}}(\mathbf{v}) = F_{{s+1}}(\mathbf{v}) \mathbf{T}_{s+1}$. Then, the fact that $\bTheta \mathbf{A} \mathbf{Q}_{(i)} = \bTheta \mathbf{A} (\mathbf{W}_{(i)} - \mathbf{Q}_{(1:i-1)} \mathbf{R}_{(1:i-1,i)}) \mathbf{R}_{(i,i)}^{-1}$ enables us to obtain the next $s-1$ columns of $\bC$ as follows:

$$ [ \bc_{s(i-2)+2}, \bc_{s(i-2)+3}, \hdots,  \bc_{s(i-2)+s} ] =  (\bTheta \bA \overline{\bW_{(i)}}  - [\bc_1, \bc_2, \hdots, \bc_{s(i-2)+1}] {\overline{\bR_{(1:i-1,i)}}) ({\overline{\overline{\bR_{(i,i)}}}}})^{-1}.$$

Here, the matrix $\overline{\bR_{(1:i-1,i)}}$ represents $\bR_{(1:i-1,i)}$ without the last column, and the matrix ${\overline{\overline{\bR_{(i,i)}}}}$ represents $\bR_{(i,i)}$ without the last row and column. Thus, the random sketching technique not only allows for the effective construction of a well-conditioned Krylov basis, but also can facilitate the determination of the Arnoldi matrix $\bH$.

\section{Numerical experiments} \label{experiments}
We test the methodology on three numerical examples:  QR factorization of a synthetically generated matrix from~\parencite[Section 5.1]{balabanov2021randomizedGS}, {the}  solution of a linear system with block GMRES, and {the} solution of an eigenvalue problem with {the} RR method. The numerical analysis of the {$s$-step RBGS-Arnoldi} algorithm and the randomized FOM method is outside the scope of this manuscript.
The RBGS process is validated by comparison with standard methods such as BCGS, BMGS and BCGS2. 
Depending on the example, finite precision arithmetic is performed in float32 or float64 format. To ensure {proper} comparison, we decided to perform the inter-block orthogonalization routines in float64, even in the unique float32 precision algorithms. Furthermore, for standard methods we perform inter-block orthogonalization with an efficient Householder QR routine, while for RBGS, this routine is combined with RCholeskyQR according to~\cref{alg:TS-CholQR}.

Several solvers are considered for the sketched least-squares problems depending on the experiment, but in all of them this step has negligible complexity and memory requirements.
Finally, in all the numerical examples, SRHT and Rademacher matrices give very similar results, so we present the results for SRHT only.

\subsection{Orthogonalization of a numerically {rank-deficient} matrix} \label{Ex1}
Take $[\bW]_{i,j} = f_{\mu_j}(x_i)$, $1\leq i \leq n$, $1 \leq j \leq m$, where
\begin{equation*}
f_\mu(x) = \frac{\sin\left (10(\mu+x) \right)}{\cos \left(100(\mu-x) \right)+1.1},
\end{equation*}
and $x_j$ and $\mu_j$ are chosen as, respectively, $n= 10^6$ and $m=300$ equally distanced points with $x_0 = \mu_0 = 0$ and $x_{10^6} = \mu_{300} = 1$. The matrix $\bW$ is partitioned into blocks of columns of size $m_p =10$. Then a block-wise QR factorization of $\bW$ is performed with {RBGS process (\cref{alg:RBGS}),} executed either in  unique float32 precision or in multi-precision, using $\bTheta$ with $k = 3000$ rows. The least-squares solver in step 2 of~\cref{alg:RBGS} is chosen either as the Householder solver, or as $20$ iterations of CG applied to the normal equation.  Thereafter, the stability behavior of RBGS is compared to {that of} the standard BCGS, BMGS, and BCGS2 algorithms executed in float32 arithmetic.

\begin{figure}[!h]
	\centering
	\begin{subfigure}{.45\textwidth}
		\centering  
		\includegraphics[width=0.8\textwidth]{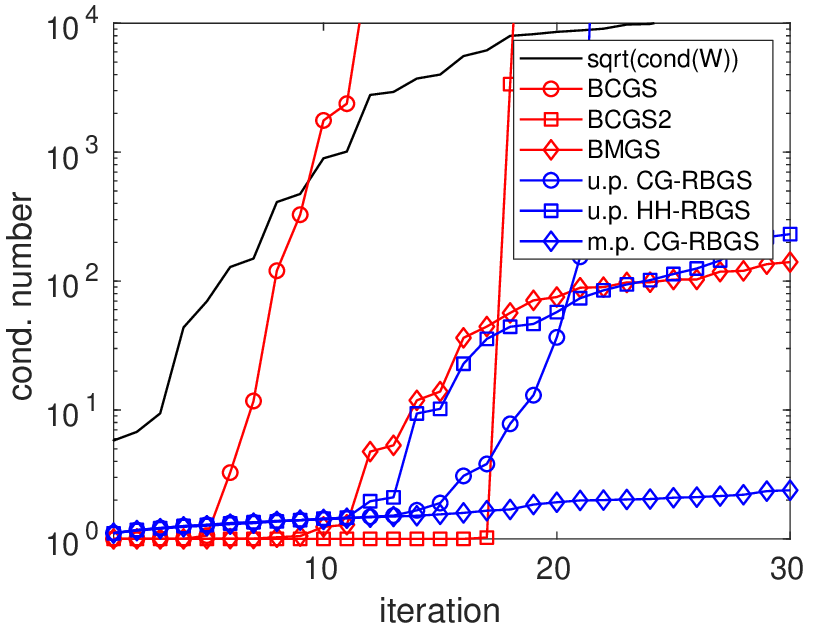}
		\caption{  Condition number of $\bQ$}
		\label{fig:Ex1_1a}
	\end{subfigure} \hspace{.03\textwidth}
	\begin{subfigure}{.45\textwidth}
		\centering
		\includegraphics[width=0.8\textwidth]{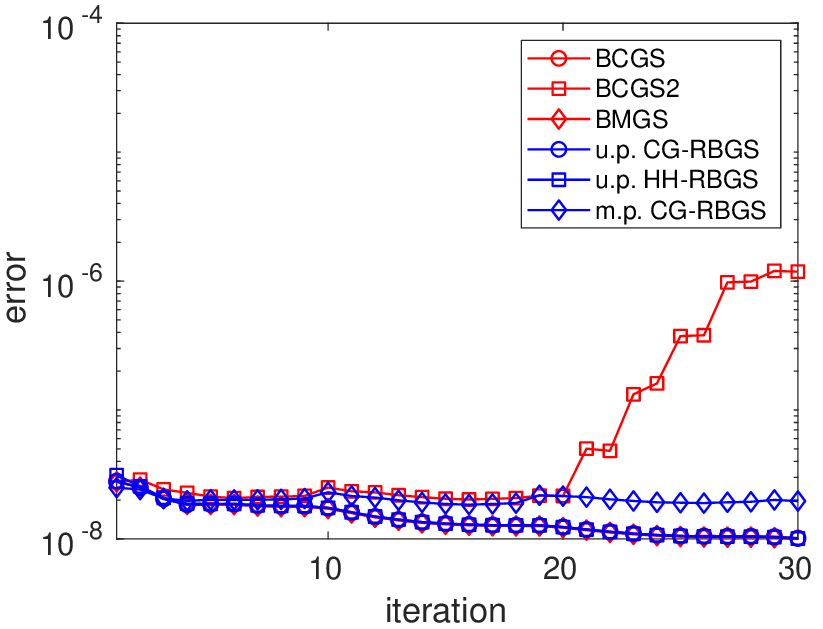}
		\caption{  Error $\|\bW - \bQ \bR\|/\|\bW\|$}
		\label{fig:Ex1_1b}
	\end{subfigure} 
	\caption{{Block QR factorization of the matrix from~\parencite[Section 5.1]{balabanov2021randomizedGS}. In the plots, {``u.p. CG-RBGS" and ``m.p. CG-RBGS" refer to the unique-precision and multi-precision RBGS algorithms, respectively, both of which perform} step 2 with $20$ iterations of CG. On the other hand, {``u.p. HH-RBGS"} refers to the unique precision RBGS {that uses a} Householder solver.}}
	\label{fig:Ex1_1}
\end{figure}

We observe a similar picture as in~\parencite[Section 5.1]{balabanov2021randomizedGS}, comparing single-vector Gram-Schmidt algorithms. According to \Cref{fig:Ex1_1a}, the BCGS and BCGS2 methods become dramatically unstable respectively  at iterations $i \geq 8$, and $i \geq 17$, with the latter iterations corresponding to $\bW_{(1:i)}$ being numerically {rank-deficient}. The BMGS algorithm exhibits better stability than the other two standard approaches. However, it still yields a Q factor with condition number of two orders of magnitude. The unique precision RBGS algorithm, using Householder least-squares solver in step 2, has {a} similar stability profile as {the} BMGS algorithm. On the other hand, the usage of $20$ iterations of CG does not provide sufficient accuracy in step 2 and, as a consequence, reduces {the} stability of the unique precision RBGS. 
In contrast to all tested methods, the multi-precision RBGS algorithm remains perfectly stable, even at iterations where $\bW_{(1:i)}$ is numerically {rank-deficient}, and outputs a Q factor with condition number close to $\mathcal{O}(1)$. Unlike in the unique precision algorithm, here the CG solver in step 2 turns out to be sufficiently accurate. Moreover,~\Cref{fig:Ex1_1b} examines the behavior of the approximation error $\|\bW_{(1:i)} - \bQ_{(1:i)} \bR_{(1:i,1:i)}\|/\|\bW_{(1:i)}\|$. We see that for all tested algorithms, besides  BCGS2, this error remains close to machine precision at all iterations.

\subsection{Solution of a linear system with block GMRES} \label{expgmres}
Consider {the} linear system 
\begin{equation}  \label{eq:ex3sys}
(\bA_{Ga}+  \alpha \bI) \bX_{Ga} = \bB,
\end{equation}
where the matrix $\bA_{Ga}$ is taken as the ``Ga41As41H72'' matrix of dimension $n = 268096$ from the  SuiteSparse  matrix  collection, $\bI$ is the $n\times n$ identity matrix, and $\alpha = 0.2$ {introduced to improve conditioning of $\bA_{Ga}$.}   The right-hand-side matrix $\bB$ is taken as $n \times m_p$ random matrix with entries being i.i.d normal random variables, and  $m_p  = 100$. Solving such {a} system could be part of the inverse subspace iteration for computing {the} eigenvalues of $\bA_{Ga}$ of the smallest magnitude. The shifted ``Ga41As41H72'' matrix has many clustered, possibly negative, eigenvalues that can bring stability issues to the solvers.

Furthermore, the system~\cref{eq:ex3sys} is preconditioned from the right by the incomplete LU factorization $\bP_{Ga}$ of $\bA_{Ga}+  \alpha \bI$ with zero level of fill-in and symmetric reverse Cuthill-McKee reordering. With this preconditioner, the final system of equations has the form $  \bA \bX = \bB, $ 
where $\bA = (\bA_{Ga}+  \alpha \bI) \bP_{Ga}$ and $\bX_{Ga} = \bP_{Ga} \bX$. {Here, {the} matrix $\bA$ is not computed explicitly but {is} considered as an {implicit} map outputting {a} product with vectors}.  This system is approximately solved with {the} GMRES method based on different versions of BGS process. We restart the GMRES method every $30$ iterations, i.e., when the dimension of the Krylov space becomes $m = 3100$. 

Here we examine the behavior of BCGS, BCGS2, BMGS and the unique precision RBGS under float32 arithmetic. There is no need to test the multi-precision RBGS as the unique precision algorithms are already able to provide {a} nearly optimal solution. The products with matrix $\bA$ and solutions of reduced GMRES least-squares problems~\cref{eq:GMRESsol} are computed in float64 format. The BGS iterations and other operations are performed in float32 format. In RBGS, the sketched least-squares problems are solved with $5$ Richardson iterations, as explained in~\cref{lssol}.

\Cref{fig:Ex3_1} provides {the} convergence of the residual error $\max_{j = 1, \hdots, 100}\|\bA \bu_j - \bb_j \|/\|\bb_j\|$. The condition number of the computed Krylov basis $\bQ_{(1:i-1)}$  at each iteration $i$ is given in~\cref{fig:Ex3_2}. We see that already starting from the first iterations, the BCGS algorithm entails dramatic instability and stagnation of the residual error. {In contrast}, the BCGS2 algorithm remains stable at all iterations, providing a Q factor that is orthonormal up to machine precision. The BMGS algorithm {shows} partial instability, resulting in deteriorated convergence of the error.  Finally, the RBGS algorithm is as stable as the BCGS2 algorithm. It provides a well-conditioned Q factor and optimal convergence of the residual error for all tested sizes of the sketching matrix. {At the same time, RBGS requires only half the flops and $\mathcal{O}(p)$ fewer global reductions than BMGS, and a quarter of the flops and at most half the global reductions of BCGS2.}

\begin{figure}[!h]
	\centering
	\begin{subfigure}{.45\textwidth}
		\centering  
		\includegraphics[width=0.8\textwidth]{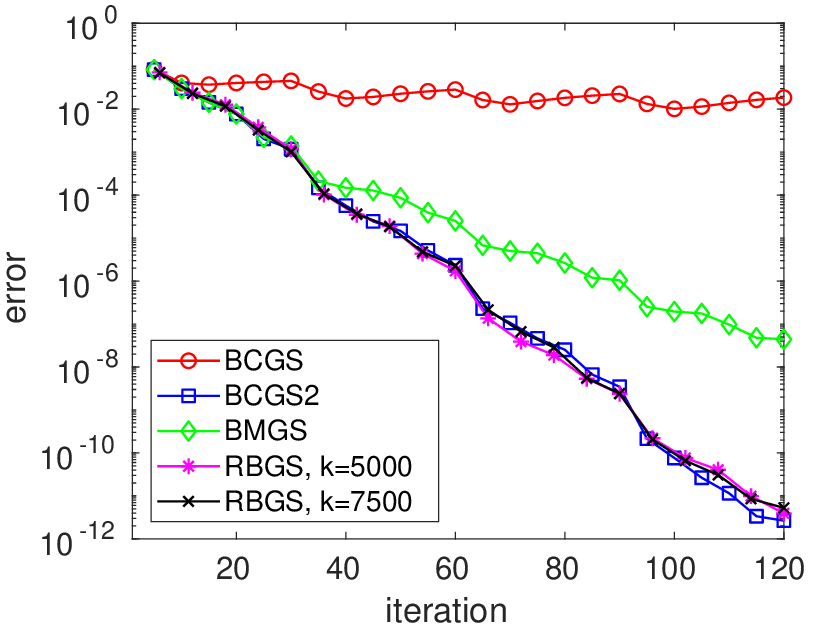}
		\caption{ Res. error $\max_{j}\|\bA \bu_j - \bb_j \|/\|\bb_j\|$}
		\label{fig:Ex2_1}
	\end{subfigure} \hspace{.03\textwidth}
	\begin{subfigure}{.45\textwidth}
		\centering
		\includegraphics[width=0.8\textwidth]{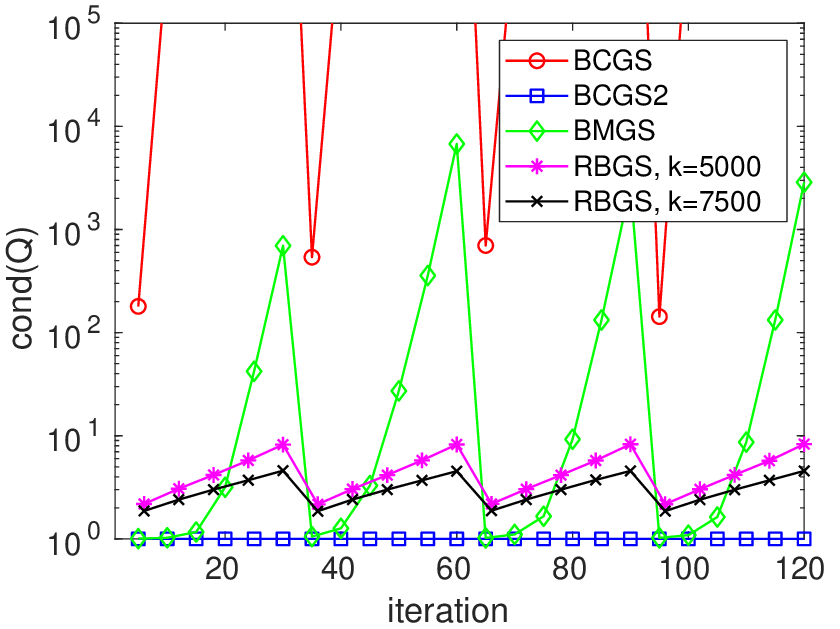}
		\caption{ Condition number of $\bQ$}
		\label{fig:Ex2_2}
	\end{subfigure}
	\caption{Solution of a linear system with GMRES.}
	\label{fig:Ex2}
\end{figure}

\subsection{Solution of an eigenvalue problem with randomized RR}
In this test, we seek the smallest (negative) eigenpairs of the ``Ga41As41H72'' matrix $\bA_{Ga}$ from~\cref{expgmres}. For this we first transform the eigenproblem to the following positive-definite one: 
\begin{equation} 
(\alpha\bI - \bA_{Ga}) \bx = \lambda \bx,
\end{equation}
where $\alpha = 1500 \approx \|\bA_{Ga}\|$, and seek the dominant eigenpairs of $\bA = \alpha\bI - \bA_{Ga}$. Then one can use the fact that $\bA_{Ga}$ and $\bA$  have {the} same eigenvectors and the associated  eigenvalues {are} related as $\lambda_{Ga} =  \alpha - \lambda$. In this case the first $50$ eigenvalues of $\bA$ are clustered inside $[1500.9~1501.3]$.  To compute the dominant eigenpairs of $\bA$ we either apply a subspace iteration or the RR method based on different BGS algorithms. The initial guess matrix $\bB$ is taken as an $n \times m_p$ Gaussian matrix of dimension $m_p= 50$.  The Arnoldi algorithm is restarted every $50$ iterations (with columns of $\bB$ chosen as the dominant Ritz vectors) which corresponds to the Krylov space of dimension $m = 2550$.  In this numerical example all arithmetic operations are performed in float64. The step 2 of {the} RBGS algorithm is again performed with $5$ Richardson iterations, as in the previous example.

\begin{figure}[!h]
	\centering
	\begin{subfigure}{.45\textwidth}
		\centering  
		\includegraphics[width=0.8\textwidth]{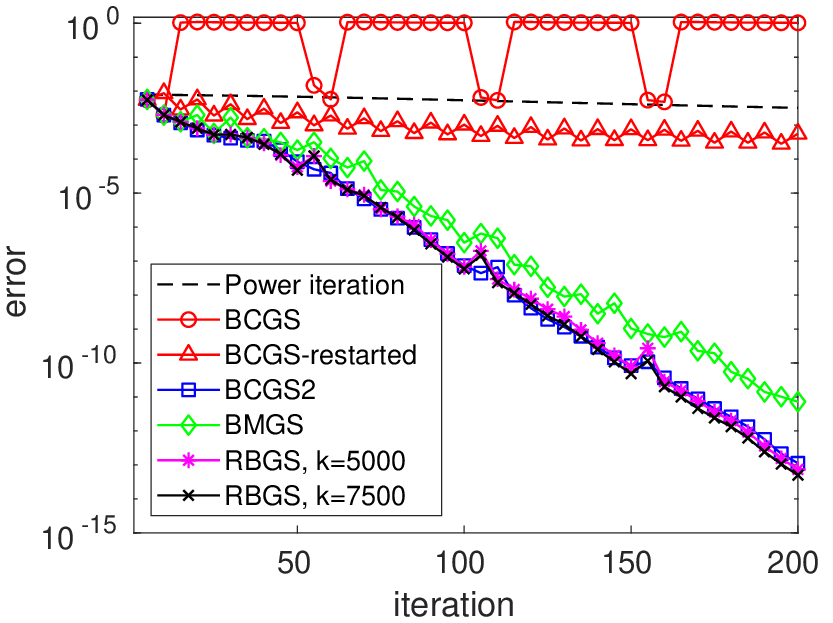}
		\caption{ Res. error $\max_{j}\|\bA \bu_j - \mu_j \bu_j \|/\|\mu_j \bu_j\|$}
		\label{fig:Ex3_1}
	\end{subfigure} \hspace{.03\textwidth}
	\begin{subfigure}{.45\textwidth}
		\centering
		\includegraphics[width=0.8\textwidth]{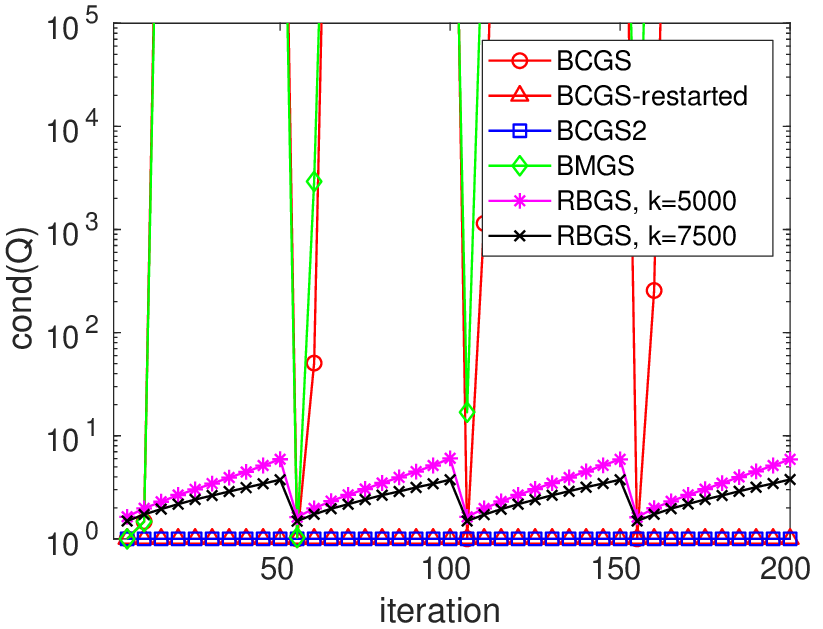}
		\caption{ Condition number of $\bQ$}
		\label{fig:Ex3_2}
	\end{subfigure}
	\caption{Solution of an eigenvalue problem with {the} RR method. In the plots, ``BCGS-restarted'' refers to {the} BCGS-Arnoldi algorithm that restarts every $5$ iterations.}
	\label{fig:Ex3}
\end{figure}

We measure the approximation error by the maximum relative residual of the first $80\%$  of computed eigenpairs (i.e., $40$ eigenpairs out of $50$).  
The convergence of the approximation error $\max_{j = 1, \hdots, 40}\|\bA \bu_j - \mu_j \bu_j \|/\|\mu_j \bu_j\|$ is depicted in~\cref{fig:Ex3_1}. The condition number of the computed Krylov basis is depicted in \cref{fig:Ex3_2}. It is revealed that the subspace iteration method yields early stagnation of the error and is unsuitable for this numerical example. On the other hand, the RR method, if stable, provides convergence of the error to machine precision. 

It is revealed that the BCGS-based Arnoldi method becomes unstable already at early iterations and produces an ill-conditioned Krylov basis with condition number close to $\mathcal{O}(u^{-1})$. Therefore this method needs to be restarted, say every $5$ iterations.  \cref{fig:Ex3_1} plots the error also for this case. It is seen that early restarting does not help as it causes a dramatic effect on the convergence of the error.  In contrast to BCGS, the BCGS2 and RBGS algorithms show perfect stability and yield an approximation that converges to machine precision in $200$ iterations. The BMGS method is unstable in the sense that it outputs an ill-conditioned Krylov basis. Despite this, it yields an approximation error that converges to machine precision, though with somewhat deteriorated rates compared to BCGS2 and RBGS.

\section{Proofs of propositions and theorems} \label{proofs}

This section contains the proofs for the stability guarantees from~\cref{stability}.
\begin{proof} [Proof of~\Cref{thm:maintheorem1}]
	We have,
	\small
	\begin{equation} \label{eq:main10}
	\| \bTheta \bhQ - \bhS\|_\Frob  \leq 1.02 u_{fine} n \| \bTheta \|_\Frob \| \bhQ\|_\Frob \leq 1.02 \sqrt{1+\varepsilon} u_{fine} n^{3/2} \| \bhQ \|_\Frob \leq 0.02 u_{crs}\| \bhQ\|:=F_1,
	\end{equation}
	\normalsize
	which implies that	
	\small		
	$$ \sigma_{min}(\bhS) - F_1 \leq \sigma_{min}(\bTheta\bhQ) \leq \sigma_{max}(\bTheta\bhQ) \leq \sigma_{max}(\bhS) + F_1.$$
	\normalsize
	By noticing that
	\small
	$$ 1-\Delta^{(p)} \leq \sqrt{1-\Delta^{(p)}} \leq  \sigma_{min}(\bhS)  \leq \sigma_{max}(\bhS) \leq \sqrt{1+\Delta^{(p)}} \leq 1+\Delta^{(p)} ,$$
	\normalsize
	and~\cref{prop:Cn2sqrtn} we deduce that $\|\bhQ \| \leq 1.57$ and that
	\small
	$$ (1+\varepsilon)^{-1/2}(1-\Delta^{(p)}- F_1) \leq \sigma_{min}(\bhQ)  \leq \sigma_{max}(\bhQ) \leq (1-\varepsilon)^{-1/2}(1+\Delta^{(p)}+ F_1).$$ \normalsize
	
	Let us now prove the second statement of the theorem. The statement clearly holds for $p = 1$. For $p\geq 2$ notice that
	\small
	\begin{equation}	\label{eq:main11}
	\begin{split}
	\| \bhR \|_\Frob 
	&\leq \sigma_{min}(\bhS)^{-1} \| \bhS \bhR\|_\Frob \leq 1.12 \| \bhS \bhR\|_\Frob \leq 1.12 (\|\bhP \|_\Frob +  \|\bhP - \bhS \bhR\|_\Frob ) \\ 
	&\leq  1.12 (1+\tilde{\Delta}^{(p)})\|\bhP\|_\Frob \leq 1.4 (1+\tilde{\Delta}^{(p)})\|\bhW\|_\Frob  \leq 1.6 \|\bhW\|_\Frob,
	\end{split}
	\end{equation}
	\normalsize
	where we used the fact that $\| \bhP \|_\Frob \leq 1.02 \|\bTheta \bhW \|_\Frob \leq {1.25} \|\bhW\|_\Frob$.
	We also have, for $1 \leq i \leq p$,
	$$ \bhW_{(i)} =  \bhQ_{(1:i)} \bhR_{(1:i,i)}  + \bDelta \bW_{(i)}, $$
	with
	\small
	\begin{align*}
	\| \bDelta \bW_{(i)}\|_\Frob &\leq \|\bhQ'_{(i)} - \bhQ_{(i)} \bhR_{(i,i)} \|_\Frob +  1.02  u_{crs} (\|\bhW_{(i)}\|_\Frob+i m_p\|\bhQ_{(1:i-1)}\|_\Frob \| \bhR_{(1:i-1,i)} \|_\Frob)\\
	&\leq 0.1 u_{crs} \|\bhQ_{(i)}\|_\Frob \| \bhR_{(i,i)} \|_\Frob +  1.02  u_{crs} (\|\bhW_{(i)}\|_\Frob+1.57 i^{3/2}m_p^{3/2}\|\bhR_{(1:i-1,i)} \|_\Frob )\\
	& \leq  1.02 u_{crs}(\|\bhW_{(i)}\|_\Frob+0.16m_p^{1/2}\| \bhR_{(i,i)}\|_\Frob +1.57 i^{3/2}m_p^{3/2} \|\bhR_{(1:i-1,i)} \|_\Frob) \\
	& \leq  1.02 u_{crs}(1.26m_p^{1/2}\|\bhW\|_\Frob +1.57 i^{3/2}m_p^{3/2} \|\bhR_{(1:i-1,i)} \|_\Frob).
	\end{align*}
	\normalsize			
	Consequently,
	\small 
	\begin{align*}
	\|\bDelta \bW\|^2_\Frob &\leq  1.02^2 u^2_{crs} \sum_{1 \leq i \leq p}  2(1.26^2m_p\|\bhW\|^2_\Frob+ 1.57^2 i^{3}m_p^{3} \|\bhR_{(1:i-1,i)} \|^2_\Frob)\\
	&\leq  1.02^2 u^2_{crs}  2(1.26^2m\|\bhW\|_\Frob^2+ 1.57^2 m^3 \| \bhR \|^2_\Frob)\\& \leq 1.02^2 u^2_{crs}2 (1.26^2 m\|\bhW\|_\Frob^2+ 2.52^2 m^3 \| \bhW \|^2_\Frob) 
	\leq \left ( 4 u_{crs} m^{3/2} \|\bhW\|_\Frob \right)^2.
	\end{align*}
	\normalsize
\end{proof}

\begin{proof} [Proof of~\Cref{thm:maintheorem1+}]
	We have, 
	\small
	$$ \bhW_{(i)} =  \bhQ_{(1:i)} \bhR_{(1:i,i)}  + \bDelta \bW_{(i)}, $$
	\normalsize
	with
	\small
	\begin{align*}
	\| \bTheta\bDelta \bW_{(i)}\|_\Frob &\leq \|\bTheta(\bhQ'_{(i)} - \bhQ_{(i)} \bhR_{(i,i)}) \|_\Frob +  1.02  \sqrt{1+\varepsilon} u_{crs} (\|\bhW_{(i)}\|_\Frob+im_p\|\bhQ_{(1:i-1)}\|_\Frob \| \bhR_{(1:i-1,i)} \|_\Frob)\\
	& \leq  1.02 u_{crs}(1.25\|\bhW_{(i)}\|_\Frob+0.16m_p^{1/2}\| \bhR_{(i,i)}\|_\Frob+2 i^{3/2}m_p^{3/2} \|\bhR_{(1:i-1,i)} \|_\Frob)\\
	& \leq  1.02 u_{crs}(1.51m_p^{1/2}  \|\bhW\|_\Frob+2 i^{3/2}m_p^{3/2} \|\bhR_{(1:i-1,i)} \|_\Frob).
	\end{align*}
	\normalsize	
	Consequently, 
	\small
	\begin{align*}
	\|\bTheta \bDelta \bW\|^2_\Frob &\leq  1.02^2 u^2_{crs}2 (1.51^2m\|\bhW\|_\Frob^2+ 2^2 m^3 \| \bhR \|^2_\Frob)\\ &\leq 1.02^2 u^2_{crs}2 (1.51^2m\|\bhW\|_\Frob^2+ 3.2^2 m^3 \| \bhW \|^2_\Frob) \leq \left ( 5 u_{crs} m^{3/2} \|\bhW\|_\Frob \right)^2.
	\end{align*}
	\normalsize
\end{proof}

\begin{proof} [Proof of~\Cref{thm:maintheorem2}]
	The proof is done by induction on $p$. Assume that~\cref{thm:maintheorem2} holds for $p = i-1 \geq 1$. Below, we show that the statement of the theorem then also holds for $p=i$.	Clearly,
	\small 
	\begin{subequations} \label{eq:mainproof1}
		\begin{align}
		\|\bhS_{(1:i-1)}\|_\Frob  &\leq \sqrt{(i-1)m_p} \sqrt{1+\Delta^{(i-1)}} \leq  1.01 \sqrt{(i-1)m_p} \\
		\| \bhP_{(i)} \| &\leq 1.02 \|\bTheta \bhW_{(i)} \| \leq {1.25} \|\bhW_{(i)}\|.
		\end{align}
	\end{subequations}
	\normalsize 
	Moreover, we have for $1\leq j \leq i$,
	\small 
	$$\| \bTheta \bhQ_{(j)} - \bhS_{(j)}\|_\Frob  \leq  1.02 u_{fine} \sqrt{1+\varepsilon} n^{3/2} \| \bhQ_{(j)} \|_\Frob \leq 0.02 u_{crs}\| \bhQ_{(j)}\|.$$
	\normalsize 
	{From the fact that $\sigma_{min}(\bhS_{(1:i-1)}) \geq \sqrt{1-\Delta^{(i-1)}} \geq 0.989 $,~\cref{thm:asssolver,eq:mainproof1}, notice that} for each $(i-1) m_p < j \leq i m_p$, the matrix $\bhS_{(1:i-1)}+ \bDelta\bS$ associated with the $j$-th column of $\bhP$ satisfies
	\small 
	\begin{equation} \label{eq:mainproof2}
	\sigma_{min}(\bhS_{(1:i-1)}+ \bDelta\bS) \geq  \sigma_{min}(\bhS_{(1:i-1)})- \| \bDelta \bS\| \geq 0.98. 
	\end{equation}
	\normalsize
	Consequently, for each $(i-1) m_p < j \leq i m_p$, the $j$-th column of $\bhR_{(1:i-1,1:i)}$ is bounded by
	$\|\bhp_{j} + \bDelta \bp \| /\sigma_{min}(\bhS_{(1:i-1)}+\bDelta \bS) \leq 1.4 \| \bhW_{j} \|$ implying that 
	\small	
	\vspace*{-1em}
	\begin{equation} \label{eq:mainproof3}
	\| \bhR_{(1:i-1,i)}\|_\Frob \leq 1.4 \| \bhW_{(i)} \|_\Frob.
	\end{equation} 	
	\normalsize
	{This fact combined with the fact (which holds by the induction hypothesis: $\Delta^{(i-1)} \leq 0.02$ and~\cref{thm:maintheorem1}) that} 
	\begin{equation} \label{eq:mainproof4}
	\|\bhQ_{(1:i-1)}\| \leq 1.5 \text{ and } \|\bhQ_{(1:i-1)}\|_\Frob \leq 1.5 (i-1)^{1/2} m_p^{1/2}, 
	\end{equation}
	and~\cref{eq:numepsembedding2}, leads to the following result:	
	\small
	\begin{gather}  
	\begin{split}
	\|\bDelta \bQ'_{(i)}\|_\Frob& = \|\bhQ'_{(i)}-(\bhW_{(i)} - \bhQ_{(1:i-1)}  \bhR_{(1:i-1,i)}) \|_\Frob \leq 1.02 u_{crs}  \| |\bhW_{(i)}|+ i m_p | \bhQ_{(1:i-1)}| |\bhR_{(1:i-1,i)} | \|_\Frob \\
	&\leq 1.02 u_{crs}  (\|\bhW_{(i)}\|_\Frob+ i m_p \| \bhQ_{(1:i-1)}\|_\Frob \|\bhR_{(1:i-1,i)} \|_\Frob)  \leq 3.2 u_{crs} i^{3/2} m_p^{3/2} \|\bhW_{(i)} \|_\Frob, \text{ \normalsize and \small } \label{eq:qprime} 
	\end{split}\\ 	
	\| \bhQ'_{(i)} \|_\Frob \leq \|\bhW_{(i)}\|_\Frob + \|\bhQ_{(1:i-1)}\| \|\bhR_{(1:i-1,i)}\|_\Frob+ \|\bDelta \bQ'_{(i)} \|_\Frob \leq 3.2   \|\bhW_{(i)} \|_\Frob. \label{eq:qprime2}
	\end{gather}
	\normalsize
	Then by~\cref{eq:numepsembedding3},
	\small
	\begin{gather}  
	\begin{split}
	\|\bTheta \bDelta \bQ'_{(i)} \|_\Frob &\leq 3.2  \sqrt{1+\varepsilon} u_{crs} i^{3/2} m_p^{3/2} \|\bhW_{(i)} \|_\Frob \leq 4 u_{crs} i^{3/2} m_p^{3/2} \|\bhW_{(i)} \|_\Frob\\
	\|\bTheta \bhQ'_{(i)}\|_\Frob &\leq \| \bTheta \bhW^{(i)}\|_\Frob +  \|\bTheta \bhQ_{(1:i-1)}\| \|\bhR_{(1:i-1,i)}\|_\Frob +  \|\bTheta \bDelta \bQ'_{{(i)}}\|_\Frob \leq 3.2 \|\bhW_{(i)}\|_\Frob.  \label{eq:Thetaqprime2}
	\end{split}
	\end{gather}  
	\normalsize			
	Denote residual matrix $\bhQ_{(1:i-1)} \bhR_{(1:i-1,1:i-1)} - \bhW_{(1:i-1)}$ by $\bhB^{(i-1)}$. Then, we have
	\small 
	\begin{gather} 
	\begin{split}
	\|\bhB^{(i-1)}\| &\leq 4 u_{crs} i^{3/2} m_p^{3/2} \|\bhW_{(1:i-1)}\|_\Frob \leq 0.01 \sigma_{min} (\bhW), \text{ and}\\
	\|\bTheta \bhB^{(i-1)}\| &\leq 5 u_{crs}i^{3/2} m_p^{3/2}  \|\bhW_{(1:i-1)}\|_\Frob \leq 0.01 \sigma_{min} (\bhW). 
	\end{split}
	\end{gather}
	\normalsize
	Now we are all set to derive the $\varepsilon'$-embedding property of $\bTheta$ for $\bhQ'_{(i)}$, which is needed  to characterize the quality of the inter-block orthogonalization step. 
	Let us first notice that, 
	\small
	\begin{gather}
	\begin{split}  
	\sigma_{min}(\bhQ'_{(i)}) &= \sigma_{min}(\bhW_{(i)} - \bhQ_{(1:i-1)}\bhR_{(1:i-1,i)} + \bDelta \bhQ'_{{(i)}}  ) \\
	&\geq \sigma_{min}(\bhW_{(i)} - \bhQ_{(1:i-1)}\bhR_{(1:i-1,i)}) - \|\bDelta \bhQ'_{{(i)}} \|  \\ 
	&\geq \sigma_{min}(\bhW_{(i)} - (\bhW_{(1:i-1)}+\bhB^{(i-1)})\bhR_{(1:i-1,1:i-1)}^{-1} \bhR_{(1:i-1,i)}) - \|\bDelta \bhQ'_{{(i)}}\|   \\ 
	& = \sigma_{min} \left ([\bhW_{(1:i-1)}+\bhB^{(i-1)},\bhW_{(i)}] \left [\begin{matrix} -\bhR_{(1:i-1,1:i-1)}^{-1} \bhR_{(1:i-1,i)} \\ \bI  \end{matrix}  \right ] \right ) - \|\bDelta \bhQ'_{{(i)}}\|   \\ 
	& \geq \sigma_{min}([\bhW_{(1:i-1)}+\bhB^{(i-1)},\bhW_{(i)}]) - \|\bDelta \bhQ'_{{(i)}}\|\\
	& \geq \sigma_{min}(\bhW) - \|\bhB^{(i-1)}\| - \|\bDelta \bhQ'_{{(i)}}\| \\
	& \geq 0.98 \sigma_{min}(\bhW).
	\end{split}
	\end{gather}
	\normalsize
	Consequently, 
	$\cond(\bhQ'_{(i)}) \leq 5 m_p^{1/2} \cond(\bhW)$.
	For any $\ba \in \mathbb{R}^{m_p}$, it holds
	\small
	\begin{align*}
	&\|\bTheta \bhQ'_{(i)} \ba \|  = \|\bTheta (\bhW_{(i)}-\bhQ_{(1:i-1)}  \bhR_{(1:i-1,i)}) \ba \| \pm \| \bTheta\bDelta \bhQ'_{{(i)}} \ba \| \\
	& = \|\bTheta (\bhW_{(i)}-  \bhW_{(1:i-1)} \bhR^{-1}_{(1:i-1,1:i-1)} \bhR_{(1:i-1,i)}) \ba \| \pm \|\bTheta  \bhB^{(i-1)}\| \|\bhR^{-1}_{(1:i-1,1:i-1)} \bhR_{(1:i-1,i)} \ba\| \pm \| \bTheta\bDelta \bhQ'_{{(i)}} \ba \| \\
	& = \left(1 \pm \sigma_{min}(\bTheta \bhW)^{-1} (\|\bTheta \bhB^{(i-1)}\|+ \| \bTheta\bDelta \bhQ'_{{(i)}}  \|) \right )\|\bTheta (\bhW_{(i)}-  \bhW_{(1:i-1)} \bhR^{-1}_{(1:i-1,1:i-1)} \bhR_{(1:i-1,i)}) \ba \|  \\
	&  = (1 \pm 0.015) \|\bTheta (\bhW_{(i)}-  \bhW_{(1:i-1)} \bhR^{-1}_{(1:i-1,1:i-1)} \bhR_{(1:i-1,i)}) \ba \|  \\
	& = (1 \pm 0.015) \sqrt{1\pm \varepsilon}\| (\bhW_{(i)}-  \bhW_{(1:i-1)} \bhR^{-1}_{(1:i-1,1:i-1)} \bhR_{(1:i-1,i)}) \ba \|  \\
	&= (1 \pm 0.015) \sqrt{1\pm \varepsilon}\left(\| (\bhW_{(i)}-  \bhQ_{(1:i-1)} \bhR_{(1:i-1,i)}) \ba \| \pm 
	\|  \bhB^{(i-1)}\| \|\bhR^{-1}_{(1:i-1,1:i-1)} \bhR_{(1:i-1,i)} \ba\| \right) \\
	& = (1 \pm 0.015) (1 \pm 0.01) \sqrt{1\pm \varepsilon}\| (\bhW_{(i)}-  \bhQ_{(1:i-1)} \bhR_{(1:i-1,i)}) \ba \|  
	\\				
	& = (1 \pm 0.015) (1 \pm 0.01) \sqrt{1\pm \varepsilon} (\| \bhQ'_{(i)} \ba \| +  \|\bDelta \bhQ'_{{(i)}} \ba\|)   
	\\		
	& = (1 \pm 0.015) (1 \pm 0.01) (1\pm 0.01) \sqrt{1\pm \varepsilon}  \| \bhQ'_{(i)} \ba \|    \\	
	&	= (1 \pm 0.036)\sqrt{1\pm \varepsilon} \| \bhQ'_{(i)} \ba \|,
	\end{align*}
	\normalsize
	which, combined with the parallelogram identity, implies that $\bTheta$ is a $\varepsilon'$-embedding for $\bhQ'_{(i)}$ with $\varepsilon' \leq 1.1\varepsilon +0.1 \leq 0.65$.
	By combining this fact with~\cref{thm:Thetaasmpts}, we obtain:  
	\small 
	$$ 1 - u_{crs} m_p^{1/2} \cond(\bhW) \leq \sigma_{min}(\bTheta \bhQ_{(i)}) \leq \sigma_{max}(\bTheta \bhQ_{(i)}) \leq 1+u_{crs} m_p^{1/2} \cond(\bhW),$$
	\normalsize
	and
	$ \| \bhQ_{(i)}\| \leq 1.8$. 
	
	Since $\| \bTheta \bhQ_{(i)} - \bhS_{(i)}\|_\Frob  \leq 0.02 u_{crs}\| \bhQ_{(i)}\|$, we have
	\small
	$$ 1 - 1.05 u_{crs} m_p^{1/2} \cond(\bhW) \leq \sigma_{min}( \bhS_{(i)}) \leq \sigma_{max}( \bhS_{(i)}) \leq 1+1.05 u_{crs} m_p^{1/2} \cond(\bhW). $$
	\normalsize
	Consequently,
	\small
	$$\| \bI - \bhS_{(i)}^\mathrm{T} \bhS_{(i)} \|_\Frob \leq m_p^{1/2} \| \bI - \bhS_{(i)}^\mathrm{T} \bhS_{(i)} \| \leq m_p^{1/2} \max_{\bx} \frac{|\|\bx\|^2 - \|\bhS_{(i)} \bx\|^2|}{\|\bx\|^2} \leq  3  u_{crs} m_p \cond(\bhW).$$
	\normalsize
	We also have \small  $\|\bhR_{(i,i)}\|_\Frob \leq 1.02 \|\bTheta \bhQ'_{(i)}\|_\Frob \leq 3.3 \|\bhW_{(i)}\|_\Frob$, \normalsize and
	\small
	\begin{align*}
	\sigma_{min} (\bhR_{(i,i)}) &\geq \sigma_{min}(\bTheta \bhQ_{(i)}\bhR_{(i,i)}) /\|\bTheta \bhQ_{(i)}\| \geq 0.99 (\sigma_{min}(\bTheta \bhQ'_{(i)}) -0.1 u_{crs}\|\bhQ_{(i)}\|\|\bhR_{(i,i)}\|) \\
	&\geq 0.58 \sigma_{min}(\bhW),  
	\end{align*}
	\normalsize
	and
	\small
	\begin{align*}
	\|\bhS_{(1:i-1)}^\mathrm{T} \bhS_{(i)}\|_\Frob &= \|\bhS_{(1:i-1)}^\mathrm{T} (\bTheta \bhQ'_{(i)}+\bDelta \bhS'_{(i)})  \bhR_{(i,i)}^{-1} + \bDelta \bhS''_{(i)}\|_\Frob \\ &\leq (\|\bhS_{(1:i-1)}^\mathrm{T} \bTheta \bhQ'_{(i)}\|_\Frob +1.01 \|\bDelta \bhS'_{(i)})\|_\Frob)/0.58 \sigma_{min}(\bhW) + \| \bDelta \bhS''_{(i)}\|_\Frob,
	\end{align*}
	\normalsize
	where $\bDelta \bhS'_{(i)} = \bTheta (\bhQ_{(i)}\bhR_{(i,i)}-\bhQ'_{(i)})$ and $\bDelta \bhS''_{(i)}  =  \bhS_{(i)}-\bTheta \bhQ_{(i)}$.
	
	By noticing that,
	\small
	\begin{align*}
	\|\bDelta \bhS'_{(i)}  \|_\Frob &\leq 0.1 u_{crs} \|\bhQ_{(i)}\| \|\bhR_{(i,i)}\| \leq 0.6 u_{crs} \| \bhW_{(i)}\|_\Frob \\
	\|\bDelta \bhS''_{(i)}  \|_\Frob& \leq 0.02 u_{crs} \| \bhQ_{(i)}\| \leq 0.04 u_{crs} \\
	\begin{split}
	\|\bhS_{(1:i-1)}^\mathrm{T} \bTheta \bhQ'_{(i)}\|_\Frob &\leq \|\bhS_{(1:i-1)}^\mathrm{T} \bTheta (\bhW_{(i)} - \bhQ_{(1:i-1)} \bhR_{(1:i-1,i)})\|_\Frob +4.04 u_{crs} i^{3/2} m_p^{3/2} \|\bhW_{(i)} \|_\Frob  \\ 
	&\leq \|\bhS_{(1:i-1)}^\mathrm{T} (\bhP_{(i)} - \bhS_{(1:i-1)} \bhR_{(1:i-1,i)})\|_\Frob +4.2 u_{crs} i^{3/2} m_p^{3/2} \|\bhW_{(i)} \|_\Frob  \\
	&\leq 5 u_{crs} i^{3/2} m_p^{3/2} \|\bhW_{(i)} \|_\Frob,
	\end{split}
	\end{align*} 
	\normalsize
	we deduce that $\|\bhS_{(1:i-1)}^\mathrm{T} \bhS_{(i)}\|_\Frob \leq 10 u_{crs} i^{3/2} m_p^{2} \cond{(\bhW)}.$

	Consequently, for $i \geq 2$,
	\small
	\begin{align*}
	(\Delta^{(i)})^2 &=  \| \bI - (\bhS_{(i)})^\mathrm{T} \bhS_{(i)} \|_\Frob^2 + 2\|\bhS_{(1:i-1)}^\mathrm{T} \bhS_{(i)}\|_\Frob^2 + (\Delta^{(i-1)})^2\\
	& \leq  (3 u_{crs} m_p \cond(\bhW))^2 + 2(10 u_{crs} i^{3/2} m_p^{2} \cond{(\bhW)} )^2 + (20 u_{crs} (i-1)^{2} m_p^{2} \cond{(\bhW)} )^2 \\
	&\leq (9+200 i^3 + 400 (i-1)^4)  (m_p^{2} u_{crs} \cond{(\bhW)} )^2 \\
	& \leq (400i^4 - 1400i^3 + 2400i^2 - 1600i + 409)  (m_p^{2} u_{crs} \cond{(\bhW)} )^2 \\
	& \leq \left (20 i^2 m_p^{2} u_{crs} \cond{(\bhW)} \right )^2
	\end{align*}
	\normalsize
	and,
	\small
	\begin{align*}
	&\|\bhS_{(1:i)} \bhR_{(1:i,i)}  - \bhP_{(i)}\|_\Frob \leq \|\bTheta \bhQ'_{(i)}- (\bhP_{(i)} -\bhS_{(1:i-1)} \bhR_{(1:i-1,i)})\|_\Frob   + \|\bhS_{(i)}\bhR_{(i,i)} - \bTheta \bhQ'_{(i)} \|_\Frob\\
	&\leq \|\bTheta \bDelta \bQ'_{(i)}\|_\Frob + \|\bhP_{(i)} - \bTheta \bhW_{(i)} \|_\Frob + \|\bTheta \bhQ_{(1:i-1)} - \bhS_{(1:i-1)} \|_\Frob \|\bhR_{(1:i-1,i)}\|\\ &~~~~~~+\|\bTheta \bhQ_{(i)} - \bhS_{(i)} \|_\Frob \|\bhR_{(i,i)} \| + \|\bTheta(\bhQ'_{(i)} - \bhQ_{(i)}\bhR_{(i,i)})\|_\Frob \\
	& \leq u_{crs} (3.2\sqrt{1+\varepsilon}  i^{3/2} m_p^{3/2} + 0.02  + 0.02  \cdot  1.5 \cdot 1.4 + 0.04  \cdot 3.3 + 0.6 ) \|\bhW_{(i)} \|_\Frob \\ &\leq 4.2 i^{3/2} m_p^{3/2} \|\bhW_{(i)} \|_\Frob.
	\end{align*}
	\normalsize
	Consequently,
	\small 
	\begin{align*}
	(\Delta^{(i)})^2 &=  \|\bhS_{(1:i)} \bhR_{(1:i,1:i)}  - \bhP_{(1:i)}  \|^2_\Frob/\|\bhP_{(1:i)}  \|^2_\Frob \\&= (\|\bhS_{(1:i)} \bhR_{(1:i,i)}  - \bhP_{(i)}  \|_\Frob^2+ \|\bhS_{(1:i-1)} \bhR_{(1:i-1,1:i-1)}  - \bhP_{(1:i-1)}  \|_\Frob^2)/\|\bhP_{(1:i)}  \|^2_\Frob  \\
	&\leq ((4.2 i^{3/2} m_p^{3/2} \|\bhW_{(i)}\|_\Frob)^2+(4.2 i^{3/2} m_p^{3/2} \|\bhW_{(1:i-1)}\|_\Frob)^2 )/\|\bhP_{(1:i)}  \|^2_\Frob \\&= (4.2 i^{3/2} m_p^{3/2} \|\bhW_{(1:i)}\|_\Frob)^2/\|\bhP_{(1:i)}\|^2_\Frob.
	\end{align*}
	\normalsize
	The proof is finished by noting that \small $\|\bhP_{(1:i)} \|_\Frob \geq (\sqrt{1-\varepsilon}-0.001)\|\bhW_{(1:i)} \|_\Frob \geq 0.7\|\bhW_{(1:i)} \|_\Frob$ \normalsize.						
\end{proof}

\begin{proof}[Proof of~\Cref{thm:aprioribound}]
	We proceed with an induction on $p$. Clearly, the statement of the proposition holds for $p=1$. Assume that $\bTheta$ is {an} $\varepsilon'$-embedding for $\bhQ_{(1:p)}$, $p = i-1 \geq 1$. This condition is sufficient for the following results in~\cref{thm:maintheorem2} and its proof to hold for $p = i$:
	\small 
	\begin{align} \label{eq:theorem2result}
	\|\bhP_{(1:i)} - \bhS_{(1:i)} \bhR_{(1:i,1:i)} \|_\Frob 
	&\leq 4.2 u_{crs} m_p^{3/2} i^{3/2} \|\bhW_{(1:i)} \|_\Frob, \\
	\Delta^{(i)} = \|\bI - \bhS_{(1:i)}^\mathrm{T} \bhS_{(1:i)}\|_\Frob 
	&\leq 20  u_{crs} m_p^2 i^2  \cond{(\bhW_{(1:i)})} \leq 0.02.
	\end{align}	
	\normalsize
	In addition, we have \small $\| \bhR_{(1:i,1:i)} \|_\Frob \leq 1.2 \|\bhW_{(1:i)} \|_\Frob$, $\|\bhQ_{(1:i-1)} \|\leq 1.9,\|\bhQ_{(i)}\| \leq 1.3 $, \normalsize and
	\small
	\begin{equation} \label{eq:aprioribound1}
	\begin{split} 
	&\|\bhW_{(1:i)} - \bhQ_{(1:i)} \bhR_{(1:i,1:i)} \|_\Frob \leq \|  \bDelta \bQ'_{(1:i)}\|_\Frob +   (\sum_{1\leq j\leq i} \|\bhQ'_{(j)} - \bhQ_{(j)} \bhR_{(j,j)} \|^2_\Frob  )^{1/2} \\
	&\leq \| \bDelta  \bQ'_{(1:i)}\|_\Frob + 0.1 u_{crs} (\sum_{1\leq j\leq i} \|\bhQ_{(j)}\|^2 \| \bhR_{(j,j)} \|^2 )^{1/2}  \\ &\leq 
	\| \bDelta  \bQ'_{(1:i)}\|_\Frob + 0.1 u_{crs} (\|\bhQ_{(i)}\| + \|\bhQ_{(1:i-1)}\|_\Frob)  \| \bhR_{(1:i,1:i)} \|_\Frob \\
	&\leq	{3.5} u_{crs} m_p^{3/2} i^{3/2} \|\bhW_{(1:i)} \|_\Frob,
	\end{split}
	\end{equation}
	\normalsize
	which can be proven directly like~\cref{eq:main11,eq:qprime}  in the proofs  of~\cref{thm:maintheorem1,thm:maintheorem2}. 
	Finally, the following holds 
	\small 
	\begin{equation} \label{eq:Thetaqs}
	\begin{split}
	\|\bhS_{(1:i)} - \bTheta \bhQ_{(1:i)}\|_\Frob \leq 1.02 u_{fine} \sqrt{1+\varepsilon} n^{3/2} \| \bhQ_{(1:i)} \|_\Frob \leq u_{crs}   0.02 \|\bhQ_{(1:i)}\|, \\
	\|\bhP_{(1:i)} - \bTheta \bhW_{(1:i)} \|_\Frob \leq 1.02 u_{fine} \sqrt{1+\varepsilon} n^{3/2} \| \bhW_{(1:i)} \|_\Frob \leq u_{crs} 0.02 \| \bhW_{(1:i)} \|.
	\end{split}
	\end{equation}
	\normalsize	
	{By using the fact that { \small $\|\bhS_{(1:i)} \| \leq 1.02$ \normalsize and \small $\sigma_{min}(\bhP_{(1:i)}) \geq  \sigma_{min}(\bTheta \bhW_{(1:i)}) - 0.02 u_{crs} \|\bhW_{(1:i)} \|$}, \normalsize along with \cref{eq:theorem2result} and the $\varepsilon$-embedding property of $\bTheta$, we get}
	\begin{equation} \label{eq:singR}
	\sigma_{min}(\bhR_{(1:i,1:i)}) \geq \frac{1}{\|\bhS_{(1:i)}\|} (\sigma_{min} (\bhP_{(1:i)}) - \|\bhP_{(1:i)} - \bhS_{(1:i)} \bhR_{(1:i,1:i)} \|_\Frob ) \geq 0.7 \sigma_{min}(\bhW_{(1:i)}).
	\end{equation}
	Properties \cref{eq:theorem2result,eq:singR,eq:Thetaqs} imply that 
	\small
	\begin{equation*} 
	\begin{split}
	\| \bTheta \bhW_{(1:i)} \bhR^{-1}_{(1:i,1:i)} - \bhS_{(1:i)} \|_\Frob 
	&\leq \| \bhP_{(1:i)} \bhR^{-1}_{(1:i,1:i)} - \bhS_{(1:i)}\|_\Frob + \|(\bhP_{(1:i)} -\bTheta \bhW_{(1:i)}) \bhR^{-1}_{(1:i,1:i)} \|_{\Frob} \\
	&\leq (\| \bhP_{(1:i)} - \bhS_{(1:i)} \bhR_{(1:i,1:i)} \|_\Frob +  \|\bhP_{(1:i)} -\bTheta \bhW_{(1:i)}\|_\Frob)\|\bhR^{-1}_{(1:i,1:i)}\| \\
	& \leq 6.2 u_{crs} m_p^{2} i^{2}  \cond(\bhW_{(1:i)}) =: F_1.
	\end{split}
	\end{equation*} 
	\normalsize
	Furthermore,
	\small
	\begin{equation*} 
	\begin{split}
	1  - \Delta^{(i)} - F_1 \leq \sigma_{min} (\bTheta \bhW_{(1:i)} \bhR^{-1}_{(1:i,1:i)}) \leq \sigma_{max} (\bTheta \bhW_{(1:i)} \bhR^{-1}_{(1:i,1:i)}) \leq 1 + \Delta^{(i)}+ F_1.
	\end{split}
	\end{equation*}
	\normalsize 
	Due to the fact that $\bTheta$ is {an} $\varepsilon$-embedding for $\bhW_{(1:i)}$, we deduce that
	\small
	\begin{equation} \label{eq:singWR}
	(1+\varepsilon)^{-1/2} (1  - \Delta^{(i)} - F_1) \leq \sigma_{min} (\bhW_{(1:i)}  \bhR^{-1}_{(1:i,1:i)}) \leq \sigma_{max} (\bhW_{(1:i)}  \bhR^{-1}_{(1:i,1:i)}) \leq (1-\varepsilon)^{-1/2} (1  + \Delta^{(i)}  + F_1).
	\end{equation} 		
	\normalsize	
	We also have from~\cref{eq:aprioribound1,eq:singR},
	\small
	\begin{equation*} 
	\|\bhW_{(1:i)} \bhR^{-1}_{(1:i,1:i)} - \bhQ_{(1:i)} \|_\Frob \leq 	3.5 u_{crs} m_p^{3/2} i^{3/2} \|\bhW_{(1:i)} \|_\Frob \|\bhR^{-1}_{(1:i,1:i)}\|
	\leq 5 u_{crs}  i^{2} m_p^{2}  \cond(\bhW_{(1:i)}) =:F_2.
	\end{equation*}
	\normalsize
	By combining this property with~\cref{eq:singWR}, we get
	\small
	\begin{equation*} 
	(1+\varepsilon)^{-1/2} (1  - \Delta_m - F_1)-F_2 \leq \sigma_{min} (\bhQ_{(1:i)}) \leq \sigma_{max} (\bhQ_{(1:i)}) \leq (1-\varepsilon)^{-1/2} (1  + \Delta_m + F_1)+F_2.
	\end{equation*} 
	\normalsize
	We conclude that
	\small
	\begin{equation} \label{eq:singQ2}
	(1+\varepsilon)^{-1/2} (1  - F_3) \leq \sigma_{min} (\bhQ_{(1:i)}) \leq \sigma_{max} (\bhQ_{(1:i)}) \leq (1-\varepsilon)^{-1/2} (1  + F_3), 
	\end{equation}
	\normalsize
	where $F_3 := \Delta_m + F_1 + \sqrt{5/4}F_2 \leq {32} u_{crs} i^2 m_p^{2} \cond{(\bhW_{(1:i)})}$,
	which in particular implies that $\|\bhQ_{(1:i)}\| \leq 1.6$ and $\sigma_{min} (\bhQ_{(1:i)}) \geq 0.86$. Furthermore, from~\cref{eq:Thetaqs}, we get
	\small
	\begin{equation} \label{eq:singthetaQ1}
	1 - F_4 \leq \sigma_{min}(\bhS_{(1:i)}) - 0.1 u_{crs} \leq  \sigma_{min} (\bTheta\bhQ_{(1:i)}) \leq \sigma_{max} (\bTheta\bhQ_{(1:i)}) \leq \sigma_{max}(\bhS_{(1:i)}) + 0.1 u_{crs} \leq  1 + F_4, 
	\end{equation}	
	\normalsize
	where $F_4 := \Delta_m +0.1 u_{crs} \leq 20.1 u_{crs} m_p^2 i^2  \cond{(\bhW_{(1:i)})}$.
	
	From~\cref{eq:singthetaQ1,eq:singQ2} it is deduced that  for any vector $\ba \in \mathbb{R}^{m}$,
	\small 
	\begin{align*}
	~~&|  \|\bhQ_{(1:i)} \ba\|^2 - \|\bTheta\bhQ_{(1:i)} \ba\|^2| \\ 
	&\leq \| \ba \|^2 \max \left\{(1-\varepsilon)^{-1}(1+F_3)^2 - (1-F_4)^2    ,(1+F_4)^2 - (1+\varepsilon)^{-1}(1-F_3)^2  \right  \} \\
	&\leq \| \ba \|^2 \max \left\{1.45 \varepsilon+2.9 F_3 +2 F_4, 1.34\varepsilon + 2.1 F_4 + 2.1 F_3\right  \} \\
	& \leq 0.74 \varepsilon'\| \ba \|^2 \leq \varepsilon'  \| \bhQ_{(1:i)} \ba \|^2,
	\end{align*}
	\normalsize
	where $\varepsilon' = 2 \varepsilon + 180 u_{crs} m_p^2 i^2 \cond(\bhW_{(1:i)}).$ By using the parallelogram identity this relation can be brought to form~\cref{eq:isometry} where $V = \mathrm{range}(\bhQ_{(1:i)} )$ and $\varepsilon =\varepsilon'$.	
	It is deduced that $\bTheta$ is a $\varepsilon'$-embedding for $\bhQ_{(1:i)}$. 
\end{proof}

\section{Concluding remarks} \label{concl}

In this work we developed a block {generalization} of the RGS process, called RBGS, to compute {a} QR factorization of a large-scale matrix. It was shown that this algorithm inherits  {the} main properties of its single-vector analogue from~\parencite{balabanov2021randomizedGS}, and in particular that it is at least as stable as the single-vector MGS process, and requires nearly half the cost of the classical BGS process in terms of flops and data passes. At the same time, RBGS is well suited for cache-based and highly parallel computational architectures because it mainly relies on matrix-matrix BLAS3 operations, in contrast to its single vector counterpart relying on matrix-vector BLAS2 operations. Like the single-vector RGS algorithm, it can be implemented using multi-precision arithmetic allowing to perform the dominant large-scale operations in precision independent of the dimension of the problem. This can {be} especially {useful} for {simulations on low-precision arithmetic architectures.}

Different strategies for treating the inter-block orthogonalization step, as well as solution of the sketched least-squares problems, have been proposed.  Special care has been taken to ensure that the computational cost of these steps is negligible compared to other steps in the RBGS {process}.  

The  stability of our algorithms was verified in numerical experiments. In particular, it was seen that the multi-precision RBGS can provide a stable QR factorization of a numerically {rank-deficient} matrix, where even the most stable standard algorithms, including BCGS2, do not work. Furthermore, the randomized Arnoldi algorithm based on RBGS showed excellent stability in the context of GMRES and RR approximation in numerical examples where BMGS showed instability and convergence degradation. These factors indicate the robustness of the RBGS algorithm. 

Our next goal is to combine RBGS with model order reduction and compression techniques for even more, possibly asymptotic, reduction of the cost of Gram-Schmidt orthogonalization and the associated Krylov methods. Another direction is the application of random sketching to increase the robustness and/or efficiency of Krylov methods that use short recurrences such as CG, BCG, Lanczos, and others.  
Furthermore, we want to improve not only the Krylov methods but also other methods that involve an orthogonalization of the approximation basis, such as block LOBPCG. 
As for the theoretical analysis, besides the characterization of the sketched Galerkin projection in the randomized RR method, we also plan to investigate the $\varepsilon$-embedding property of the sketching matrix for the computed Krylov space in the presence of rounding errors. Despite the fact that this property has been thoroughly verified in numerical experiments, it still remains unproven for the randomized Arnoldi algorithm.

\section{Acknowledgments}
This project has received funding from the European Research Council (ERC) under the European Union's Horizon 2020 research and innovation program (grant agreement No 810367).

\printbibliography

\newpage
\section*{Supplementary materials: solving least-squares problem at step 2 of BRGS}
\renewcommand*{\thesection}{SM}
\begin{refsection}
	This section discusses iterative methods for solving almost orthogonal least-squares problem at step 2 of the BRGS algorithm.
	
	\subsection{Richardson iterations}
	\Cref{alg:Richardson} describes, perhaps, the simplest iterative method for obtaining a solution in step 2. When the least-squares problem is seen as orthogonalization {of columns of} $\bP_{(i)}$ to $\bS_{(1:i-1)}$, this algorithm can be interpreted as nothing more {than the classical Gram-Schmidt (CGS) algorithm with $l$ re-orthogonalizations.}   In more general terms, \Cref{alg:Richardson} is Richardson method applied to the normal system of equations
	\begin{equation} \label{eq:normalsys}
	\left (\bS_{(1:i-1)}^\mathrm{T} \bS_{(1:i-1)} \right ) \bX =  \bS_{(1:i-1)}^\mathrm{T} \bP_{(i)}. 
	\end{equation}   
	Assuming that $l = \mathcal{O}(1)$, it requires $\mathcal{O}(m_pmk)$ flops for computing $\bR_{(1:i-1,i)} = \bX$, which makes the cost of step 2 negligible in comparison to other costs of the RBGS algorithm and in particular to the step 3 of RBGS, which has complexity  $\mathcal{O}(m_pmn)$. Finally, we note that the RBGS algorithm based on $l=1$ Richardson iteration in step 2 is exactly equivalent to the BCGS process for the orthogonalization with respect to the sketched inner product. In general, there is a strong connection between the RBGS algorithm using $l$ Richardson iterations and BCGS using $l$ re-orthogonalizations. 
	
	\begin{algorithm}[h] \caption{Step 2 of RBGS, using CGS with $l$ re-orthogonalizations} \label{alg:Richardson}
		\begin{algorithmic}
			\STATE{\textbf{Given:} $\bP_{(i)}$, $\bS_{(1:i-1)}$} 
			\STATE{\textbf{Output}:  $\bR_{(1:i-1,i)} = \bX$}
			\STATE{1. Set $\bX = \bnull_{(1:i-1,1)}$.} 
			\FOR{$1:l$} 
			\STATE{2. $\bX \leftarrow \bX + \bS_{(1:i-1)}^\mathrm{T} \left (\bP_{(i)} -\bS_{(1:i-1)} \bX \right )$.}	
			\ENDFOR
		\end{algorithmic}
	\end{algorithm}

	Let $\bY = [\by_1, \hdots, \by_{m_p}]$ be the exact solution to~\cref{eq:normalsys} and $\bR_{(1:i-1,i)} = \bX =  [\bx_1, \hdots, \bx_{m_p}]$ be the solution computed with~\cref{alg:Richardson}.  In exact arithmetic, we have: 
	\begin{equation} \label{eq:Richardson}
	\|\bx_j - \by_j\| \leq \| \bI - \bS_{(1:i-1)}^\mathrm{T}\bS_{(1:i-1)}\|^l \|\by_j\| = (\Delta^{(i-1)})^l \|\by_j\|,~1\leq j \leq m_p,
	\end{equation}
	where {$\Delta^{(i-1)} = \| \bI - \bS_{(1:i-1)}^\mathrm{T}\bS_{(1:i-1)} \|_\Frob$ measures the orthogonality of the sketch of Q factor.}
	This result can be easily extended to finite precision arithmetic. 

	\subsection{Iterations of the modified GS process}
	Another way to compute the solution in step 2 of~\Cref{alg:RBGS} is to use the modified GS (MGS) process, applied $l$ times. This algorithm can provide an accurate solution in fewer iterations than~\cref{alg:Richardson}, though it can have a higher cost per iteration from a performance standpoint.  This drawback can be remedied by appealing to the block version of MGS depicted in~\cref{alg:lMGS}. The RBGS algorithm based on~\cref{alg:lMGS} can be linked to BMGS algorithm, in which the  $\ell_2$-inner products are replaced by $\langle \bTheta \cdot, \bTheta \cdot \rangle$. In particular, for the case $l=1$ the two algorithms are essentially equivalent. 
	The numerical analysis of~\cref{alg:lMGS} is beyond the scope of this manuscript. 
	\begin{algorithm}[h] \caption{Step 2 of RBGS, using BMGS with $l$ re-orthogonalizations}  \label{alg:lMGS}
		\begin{algorithmic}
			\STATE{\textbf{Given:} $\bP_{(i)}$, $\bS_{(1:i-1)}$} 
			\STATE{\textbf{Output}: $\bR_{(1:i-1,i)} = \bX = [\bX_{(1)}, \hdots, \bX_{(i-1)}]$}
			\STATE{1. Set $\bX = \bnull_{(1:i-1,1)}$.} 
			\FOR{$1:l$} 
			\FOR{$j = 1:i-1$} 
			\STATE{2. $\bP_{(i)} \leftarrow \bP_{(i)} - \bS_{(j)} \bX_{(j)}$.}
			\STATE{3. $\bX_{(j)} \leftarrow \bX_{(j)} + \bS_{(j)}^\mathrm{T} \bP_{(i)}$.}
			\ENDFOR
			\ENDFOR
		\end{algorithmic}
	\end{algorithm}
	
	\subsection{Other methods}
	In principle, we can compute $\bX$ from the normal equation~\cref{eq:normalsys} using any suitable iterative method such as Conjugate Gradient or GMRES. In this case, the normal matrix can be operated with as an implicit map outputting products with vectors using  $\mathcal{O}(m_pmk)$ flops,  making the cost of the iterative method similar to that of a Richardson iteration. We have $\kappa : = \cond(\bS_{(1:i-1)}^\mathrm{T} \bS_{(1:i-1)}) = \cond(\bS_{(1:i-1)})^2 \leq \left (\frac{1+\Delta^{(i-1)}}{1-\Delta^{(i-1)}} \right )^2$. Then the standard results on the convergence of Krylov methods, given, for example, in~\parencite[Section 6.11]{saad2003iterative}, guarantee that the solution after $l$ iterations of CG satisfies the relation
	$$ \| \bS_{(1:i-1)} (\bx_j - \by_j)  \| \leq 2 \left ( \frac{\sqrt{\kappa}-1}{\sqrt{\kappa}+1} \right )^l \|\by_j  \| \leq 2 (\Delta^{(i-1)})^l \|\by_j  \|.$$ 
	This result is similar to~\cref{eq:Richardson} for the Richardson iterations. In practice, however, the Krylov methods are expected to be more accurate and robust.

\section*{Supplementary materials: brief analysis of the accuracy of the randomized RR approximation} \label{skGalerkin}
This section provides a characterization of the accuracy of the randomized RR (or Galerkin) approximation from~\cref{rRR}. We proceed with reformulation of the methodology in terms of projection operators similarly to~\parencite[Section 4.3]{saad2011numerical} for classical methods. Our analysis will be based on the fact that $\bTheta$ satisfies the $\varepsilon$-embedding property for some finite collection $\mathcal{V}$ of fixed low-dimensional subspaces $V$. This property can be satisfied with probability at least $1- \delta \#\mathcal{V}$, if $\bTheta$ is an $(\varepsilon,\delta , d)$ oblivious $\ell_2$-subspace embedding, with $d = \max_{V \in \mathcal{V}} \mathrm{dim}(V)$.

Let $\mathcal{K}$ be an approximation space, which does not necessarily have to be a Krylov space. 
Let $\bPi_{\mathcal{K}}$ denote the $\ell_2$-orthogonal projector onto $\mathcal{K}$:
\begin{equation}
\forall \bw \in \mathbb{R}^n,~\bPi_{\mathcal{K}} \bw = \arg\min_{\bv \in \mathcal{K}} \| \bw - \bv \|.
\end{equation}
Then the Galerkin orthogonality condition~\cref{eq:Galerkin} (with $\mathcal{K} = Q_{p-1}$) can be expressed as 
$$ \bPi_{\mathcal{K}} \br(\bu,\mu\bu) = 0,$$
or, equivalently, 
$$ \bPi_{\mathcal{K}} \bA \bPi_{\mathcal{K}}  \bu = \mu\bu.$$
In other words, the classical RR method can be interpreted as approximation of eigenpairs of $\bA$ by the eigenpairs of \emph{approximate} operator $\bPi_{\mathcal{K}} \bA \bPi_{\mathcal{K}}$. 

Similarly, the sketched  Galerkin orthogonality condition~\cref{eq:skGalerkin} (with $\mathcal{K} = Q_{p-1}$) can be expressed as 
\begin{equation} \label{eq:skGalerkin2}
\bPi^\bTheta_{\mathcal{K}} \br(\bu,\mu\bu) = 0,
\end{equation}
or, equivalently, 
$$ \bPi^\bTheta_{\mathcal{K}} \bA \bPi^\bTheta_{\mathcal{K}}  \bu = \mu\bu,$$
where $\bPi^\bTheta_{\mathcal{K}}$ is an orthogonal projector onto $\mathcal{K}$ with respect to the sketched inner product $\langle \bTheta \cdot, \bTheta \cdot \rangle$, i.e.,
\begin{equation} 
\forall \bw \in \mathbb{R}^n,~\bPi^\bTheta_{\mathcal{K}} \bw = \arg\min_{\bv \in \mathcal{K}} \| \bTheta(\bw- \bv) \|,
\end{equation}
or, in matrix form, $$\bPi^\bTheta_{\mathcal{K}} = \bQ (\bTheta \bQ)^\dagger \bTheta,$$ where $\bQ$ is a matrix whose columns compose a basis for $\mathcal{K}$.
We see that the randomized RR method corresponds to taking the {approximate} operator as $\bPi^\bTheta_{\mathcal{K}} \bA \bPi^\bTheta_{\mathcal{K}}$ instead of $\bPi_{\mathcal{K}} \bA \bPi_{\mathcal{K}}$.
Next we shall provide an upper bound for the residual norm of the exact eigenpair with respect to such approximate operator. For this, it is necessary to first establish a characterization of the sketched orthogonal projector, given in~\cref{thm:lemma1}.
\begin{lemma} \label{thm:lemma1}
	Let $\bw \in \mathbb{R}^n$. If $\bTheta$ is $\varepsilon$-embedding for $V = \mathcal{K}+\mathrm{span}(\bw)$, then
	\begin{subequations}
		\begin{align}
		&\|(\bI -\bPi_{\mathcal{K}}) \bw\| \leq \|(\bI -\bPi^\bTheta_{\mathcal{K}}) \bw\| \leq \sqrt{\frac{1+\varepsilon}{1-\varepsilon}}\|(\bI -\bPi_{\mathcal{K}}) \bw\|, \\ 
		&\frac{1}{1+\varepsilon}( \|\bPi_{\mathcal{K}}\bw\|^2  - 2\varepsilon\| \bw\|^2 ) \leq \|\bPi^\bTheta_{\mathcal{K}} \bw\|^2  \leq \frac{1}{1-\varepsilon}( \|\bPi_{\mathcal{K}}\bw\|^2  + 2\varepsilon\| \bw\|^2 ), \label{eq:lemma1_2} \\
		\intertext{Moreover, if $\bw \in \mathcal{K}$, then}
		&~~~~~~~~~~~~~~~~~~~~~~~~~~~~{\bPi^\bTheta_{\mathcal{K}}  \bw = \bw.} \label{eq:lemma1_3}
		\end{align}
	\end{subequations}
	\begin{proof}
		By definitions of $\bTheta$, $\bPi_{\mathcal{K}}$ and $\bPi^\bTheta_{\mathcal{K}}$,   we have
		$$ \|(\bI -\bPi_{\mathcal{K}}) \bu\|^2 \leq \|(\bI -\bPi^\bTheta_{\mathcal{K}}) \bu\|^2 \leq \frac{1}{1-\varepsilon}(\|\bTheta(\bI -\bPi^\bTheta_{\mathcal{K}}) \bu\|)^2 \leq  \frac{1}{1-\varepsilon}(\|\bTheta(\bI -\bPi_{\mathcal{K}}) \bu\|)^2 \leq \frac{1+\varepsilon}{1-\varepsilon}(\|(\bI -\bPi_{\mathcal{K}}) \bu\|)^2,$$ 
		which gives the first inequality. 
		Furthermore, we also have the following relation 
		$$\| \bTheta \bPi^\bTheta_{\mathcal{K}} \bw\|^2 = \|\bTheta \bw\|^2 - \| \bTheta(\bI-\bPi^\bTheta_{\mathcal{K}}) \bw\|^2, $$
		which implies that
		$$ \|\bTheta \bPi^\bTheta_{\mathcal{K}}\bw\|^2 \leq (1+\varepsilon)\| \bw\|^2 - \| \bTheta(\bI-\bPi^\bTheta_{\mathcal{K}}) \bw\|^2 \leq (1+\varepsilon)\| \bw\|^2 - (1-\varepsilon)\|(\bI-\bPi_{\mathcal{K}}) \bw\|^2 \leq \| \bPi_{\mathcal{K}} \bw\|^2 + 2 \varepsilon \|\bw\|^2.   $$
		Similarly,
		$$ \|\bTheta \bPi^\bTheta_{\mathcal{K}}\bw\|^2 \geq (1-\varepsilon)\| \bw\|^2 - \| \bTheta(\bI-\bPi_{\mathcal{K}}) \bw\|^2 \geq (1-\varepsilon)\| \bw\|^2 - (1+\varepsilon)\|(\bI-\bPi_{\mathcal{K}}) \bw\|^2 \geq \| \bPi_{\mathcal{K}} \bw\|^2 - 2 \varepsilon \|\bw\|^2.   $$
		By using the relation, 
		$$ \frac{1}{1+\varepsilon}\|\bTheta \bPi^\bTheta_{\mathcal{K}}\bw\|^2 \leq \|\bPi^\bTheta_{\mathcal{K}}\bw\|^2 \leq \frac{1}{1-\varepsilon} \|\bTheta \bPi^\bTheta_{\mathcal{K}}\bw\|^2,   $$ 
		we obtain~\cref{eq:lemma1_2}.
		Finally, if $\bw \in \mathcal{K}$, then, by definition, $\bPi^\bTheta_{\mathcal{K}} \bw = \arg \min_{\bv \in \mathcal{K}} \|\bTheta(\bw - \bv)\| = \bw$, which finishes the proof.
	\end{proof}
\end{lemma}

\begin{theorem}(extension of~\parencite[Theorem 4.3]{saad2011numerical}) \label{thm:skGalerkin0}
	Assume that $\bTheta$ is an $\varepsilon$-embedding for $V = \mathcal{K}+\bA\mathcal{K}+\mathrm{span}(\bx)$.
	Let $$\gamma := \max_{\bv \in \mathcal{K}+\mathrm{span}(\bx), \|\bv\| =1 } \|\bTheta \bPi^\bTheta_{\mathcal{K}} \bA(\bI - \bPi^\bTheta_{\mathcal{K}}) \bv\|.$$ Then the sketched residual norms of the pairs $(\lambda, \bPi^\bTheta_{\mathcal{K}} \bx)$ and $(\lambda, \bx)$ for the linear operator $\bA_m = \bPi^\bTheta_{\mathcal{K}} \bA \bPi^\bTheta_{\mathcal{K}}$ satisfy, respectively,
	\begin{align*}
	&\| \bTheta (\bA_m - \lambda \bI)\bPi^\bTheta_{\mathcal{K}} \bx \| \leq \gamma \| (\bI -\bPi^\bTheta_{\mathcal{K}})\bx \| \\
	&\| \bTheta (\bA_m - \lambda \bI) \bx \| \leq \sqrt{\lambda^2+\gamma^2} \| (\bI -\bPi^\bTheta_{\mathcal{K}})\bx \|.
	\end{align*}
	\begin{proof}
		The proof directly follows that of~\cite[Theorem 4.3]{saad2011numerical} replacing $\ell_2$-inner products and norms by the sketched ones.
		We have,
		\begin{align*}
		\|\bTheta (\bA_m - \lambda \bI) \bPi^\bTheta_\mathcal{K} \bx \| &= \|\bTheta\bPi^\bTheta_\mathcal{K}(\bA - \lambda \bI) (\bx - (\bI-\bPi^\bTheta_\mathcal{K}) \bx \|  =\|\bTheta\bPi^\bTheta_\mathcal{K}(\bA - \lambda \bI) (\bI-\bPi^\bTheta_\mathcal{K}) \bx \| \\
		&= \|\bTheta\bPi^\bTheta_\mathcal{K}(\bA - \lambda \bI) (\bI-\bPi^\bTheta_\mathcal{K} ) (\bI-\bPi^\bTheta_\mathcal{K}) \bx \| \\
		& \leq \gamma \| (\bI - \bPi^\bTheta_\mathcal{K}) \bx \|.
		\end{align*}
		Notice that $\bA_m(\bI-\bPi^\bTheta_\mathcal{K}) \bx = \bnull$. Consequently,	
		$$ (\bA_m - \lambda \bI)  \bx  =  (\bA_m - \lambda \bI)\bPi^\bTheta_\mathcal{K} \bx + (\bA_m - \lambda \bI)  (\bx - \bPi^\bTheta_\mathcal{K} \bx)  = (\bA_m - \lambda \bI)\bPi^\bTheta_\mathcal{K} \bx + \lambda (\bI - \bPi^\bTheta_\mathcal{K}) \bx.  $$
		Furthermore, by using the fact that the two vectors on the right hand
		side are orthogonal with respect to the sketched inner product, we obtain 
		$$\| \bTheta (\bA_m - \lambda \bI)  \bx\|^2  =  \|\bTheta (\bA_m - \lambda \bI)\bPi^\bTheta_\mathcal{K} \bx\|^2 + |\lambda|^2  \|\bTheta (\bI - \bPi^\bTheta_\mathcal{K}) \bx \|^2 \leq (\gamma^2+ |\lambda|^2 ) \| (\bI - \bPi^\bTheta_\mathcal{K}) \bx \|^2.$$		
	\end{proof}
\end{theorem}
\Cref{thm:skGalerkin0} together with~\cref{thm:lemma1} imply that  
\begin{align*}
&\| (\bA_m - \lambda \bI)\bPi^\bTheta_{\mathcal{K}} \bx \| \leq (1+\mathcal{O}(\varepsilon))\gamma \| (\bI -\bPi_{\mathcal{K}})\bx \| \\
&\|  (\bA_m - \lambda \bI) \bx \| \leq (1+\mathcal{O}(\varepsilon)) \sqrt{\lambda^2+\gamma^2} \| (\bI -\bPi_{\mathcal{K}})\bx \|.
\end{align*}
Furthermore, according to~\cref{thm:lemma1}, we also have
$$ \gamma \leq \| \bPi_{\mathcal{K}} \bA(\bI - \bPi_{\mathcal{K}})\| + \mathcal{O}(\varepsilon)\| \bA\| \leq (1+\mathcal{O}(\varepsilon))\|\bA\|.$$
Thus, \cref{thm:skGalerkin0} guarantees a good approximation if the distance $\| (\bI - \bPi_ \mathcal{K}) \bx \| $ is small and if the approximate eigenproblem is well conditioned. Interestingly, by taking $\bTheta$ as an identity in \Cref{thm:skGalerkin0}, we exactly recover Theorem 4.3 in~\parencite{saad2011numerical}.

On the contrary, the following result bounds the residual norm of the approximate eigenpair $(\mu, \bu)$ in~\cref{eq:skGalerkin2} with respect to the exact operator.  

\begin{theorem} \label{thm:skGalerkin}
	Let $(\mu,\bu)$ be a solution to~\cref{eq:skGalerkin2}. If  $\bTheta$ is $\varepsilon$-embedding for $V = \mathcal{K}+\bA\mathcal{K}$, then we have
	\begin{equation*}
	\| (\bA - \mu \bI) \bu \| \leq (1+\mathcal{O}(\varepsilon))\|(\bI - \bPi_{\mathcal{K}})\bA\bPi_{\mathcal{K}} \| \|\bu\|,
	\end{equation*}
	Moreover, if $\mathcal{K}$ is a Krylov space $Q_{p-1}$, then 
	$$\| (\bA - \mu \bI) \bu \| \leq (1+\mathcal{O}(\varepsilon))\|(\bI - \bPi_{Q_{p-1}})\bA\bPi_{Q_{p-1}} \| \|(\bI -\bPi_{Q_{p-2}} )\bu\|,$$
	\begin{proof}
		We have, 
		$$ \|\bTheta \br(\bu,\mu \bu) \|^2 = \|\bTheta \bPi^\bTheta _{\mathcal{K}}\br(\bu,\mu \bu) \|^2 +  \|\bTheta (\bI-\bPi_{\mathcal{K}}^\bTheta)\br(\bu,\mu \bu) \|^2 = \|\bTheta (\bI-\bPi_{\mathcal{K}}^\bTheta)\br(\bu,\mu \bu) \|^2.$$
		Consequently, by the $\varepsilon$-embedding property of $\bTheta$, we obtain
		\begin{align*}
		\|\br(\bu,\mu \bu) \| &\leq \frac{1}{\sqrt{1-\varepsilon}} \|\bTheta \br(\bu,\mu \bu) \| =  \frac{1}{\sqrt{1-\varepsilon}} \|\bTheta (\bI-\bPi_{\mathcal{K}}^\bTheta)\br(\bu,\mu \bu) \| \leq \frac{1}{\sqrt{1-\varepsilon}} \|\bTheta (\bI-\bPi_{\mathcal{K}})\br(\bu,\mu \bu) \| \\
		& \leq \sqrt{\frac{1+\varepsilon}{{1-\varepsilon}}} \| (\bI-\bPi_{\mathcal{K}})\br(\bu,\mu \bu) \|  = \sqrt{\frac{1+\varepsilon}{{1-\varepsilon}}} \| (\bI-\bPi_{\mathcal{K}})\bA \bu\| = \sqrt{\frac{1+\varepsilon}{{1-\varepsilon}}} \| (\bI-\bPi_{\mathcal{K}})\bA \bPi_{\mathcal{K}} \bu\|.
		\end{align*}
		which gives theorem's first inequality. Moreover, notice that 
		\begin{align*}
		(\bI-\bPi_{{Q}_{p-1}})\bA \bPi_{{Q}_{p-1}} \bu = (\bI-\bPi_{{Q}_{p-1}})\bA (\bI -\bPi_{Q_{p-2}} ) \bu,
		\end{align*}
		which gives theorem's second inequality.
	\end{proof}
\end{theorem}

\begin{corollary}
	Assume that $\bTheta$ is an $\varepsilon$-embedding for $\mathcal{K}$. If $\mathcal{K}$ is invariant under $\bA$ then each pair $(\mu, \bu)$ that satisfies the sketched Galerkin orthogonality condition~\cref{eq:skGalerkin2} is an eigenpair of $\bA$.  
\end{corollary}

\Cref{thm:skGalerkin} implies two important results. The first one is that the sketched Galerkin projection provides an accurate result in terms of the residual error when the approximation space $\mathcal{K}$ approximates well the range of $\bA$ or when $\mathcal{K}$ is close to an invariant space. In particular, it follows that if we have  $\|(\bI - \bPi_{\mathcal{K}})\bA\bPi_{\mathcal{K}}\| \leq \tau |\mu_{m}|$, then the dominant $m$ computed  eigenpairs $(\mu,\bu)$ have relative error less than $\tau$:
\begin{equation} \label{eq:galerror}
\frac{\| (\bA - \mu \bI) \bu \| }{\|\mu \bu \|} \leq (1+\mathcal{O}(\varepsilon))\tau.
\end{equation}
At the same time,~\cref{thm:skGalerkin} also proves that if $\mathcal{K}$ represents a Krylov subspace $Q_{p-1}$, then the convergence of eigenvectors implies convergence of their residuals to zero. 
More rigorously, if $$ \|\bu^{(p)} - \bu^{(p-1)}\| \leq \tau |\mu_{m}| \|\bu^{(p)}\|, $$
where $\bu^{(p)}$ denotes the sketched Galerkin solution in $Q_{p-1}$, then~\cref{eq:galerror} holds.

Unfortunately, the fact that the residual norm is small does not guarantee that the eigenpair is accurate due to possibility of bad conditioning of the eigenvalue.  If $\bA$ is symmetric (or Hermitian for the complex case), this issue can be circumvented for classical Galerkin approximation,  by appealing to the Min-Max principle and the Courant characterization~\parencite{saad2011numerical}. Perhaps, similar ideas can be also used for the sketched  Galerkin projection. The sketched Rayleigh quotient $R(\bw)$ can be naturally defined as a number that minimizes $\|\bTheta(\bA \bw - R(\bw) \bw)\|$. It is easy to see that $ R(\bw) = \frac{\langle \bTheta\bw, \bTheta\bA \bw \rangle}{\langle \bTheta \bw,\bTheta \bw  \rangle}.$  Furthermore, it can be shown that the solution $(\mu,\bu)$ to~\cref{eq:skGalerkin2} and the exact solution $(\lambda,\bx)$ satisfy $\mu = R(\bu)$ and $\lambda = R(\bx)$. We can also show that $|R(\bPi_{\mathcal{K}}\bx) - \lambda|$ is small, if $\|\bI - \bPi_{\mathcal{K}} \bx\|$ is small. The main difficulty now becomes how to use the Min-Max principle, since the sketched operator $\bPi^\bTheta_{\mathcal{K}} \bA \bPi^\bTheta_{\mathcal{K}}$ can be non-symmetric, even when $\bA$ is symmetric.

\printbibliography[heading=subbibliography]
\end{refsection}

\end{document}